\documentclass[12pt,oneside,reqno]{amsart}
\usepackage{amssymb}
\usepackage{}
\usepackage{txfonts}
\usepackage{bbm}
\usepackage{amsmath}
\usepackage{graphicx}
\usepackage{mathrsfs}
\usepackage{stmaryrd}
\usepackage{amsfonts}
\usepackage{enumerate,amsmath,amssymb,amsthm}

\numberwithin{equation}{section}

\newcommand{\be}{\begin{eqnarray}}
\newcommand{\ee}{\end{eqnarray}}
\newcommand{\ce}{\begin{eqnarray*}}
\newcommand{\de}{\end{eqnarray*}}
\newtheorem{theorem}{Theorem}[section]
\newtheorem{lemma}[theorem]{Lemma}
\newtheorem{remark}[theorem]{Remark}
\newtheorem{definition}[theorem]{Definition}
\newtheorem{proposition}[theorem]{Proposition}
\newtheorem{Examples}[theorem]{Example}
\newtheorem{corollary}[theorem]{Corollary}

\def\eps{\varepsilon}

\def\p{\partial}

\def\[{{\Big[}}
\def\]{{\Big]}}
\def\<{{\langle}}
\def\>{{\rangle}}
\def\({{\Big(}}
\def\){{\Big)}}

\def\bx{{\mathbf{x}}}

\def\dif{{\mathord{{\rm d}}}}

\def\no{\nonumber}
\def\={&\!\!=\!\!&}
\def\bt{\begin{theorem}}
\def\et{\end{theorem}}
\def\bl{\begin{lemma}}
\def\el{\end{lemma}}
\def\br{\begin{remark}}
\def\er{\end{remark}}

\def\bd{\begin{definition}}
\def\ed{\end{definition}}
\def\bp{\begin{proposition}}
\def\ep{\end{proposition}}
\def\bc{\begin{corollary}}
\def\ec{\end{corollary}}
\def\bx{\begin{Examples}}
\def\ex{\end{Examples}}

\def\cB{{\mathcal B}}

\def\mH{{\mathbb H}}

\def\mK{{\mathbb K}}

\def\mN{{\mathbb N}}

\def\mR{{\mathbb R}}

\def\sB{{\mathscr B}}

\def\sK{{\mathscr K}}
\def\sL{{\mathscr L}}

\def\geq{\geqslant}
\def\leq{\leqslant}

\allowdisplaybreaks

\begin{document}
\title{Heat kernel estimates for critical fractional diffusion operator}

\date{}
\author{Longjie Xie and Xicheng Zhang}

\address{Longjie Xie: School of Mathematics and Statistics, Wuhan University,
Wuhan, Hubei 430072, P.R.China\\
Email: xlj.98@qq.com
}

\address{Xicheng Zhang:
School of Mathematics and Statistics, Wuhan University,
Wuhan, Hubei 430072, P.R.China\\
Email: XichengZhang@gmail.com
 }

\begin{abstract}
In this work we construct the heat kernel of the $\frac{1}{2}$-order Laplacian perturbed by the first-order gradient term
in H\"older space and the zero-order potential term in generalized Kato's class,
and obtain sharp two-sided estimates as well as the gradient estimate of the heat kernel.

\end{abstract}
\maketitle

\section{Introduction and Main Result}

For $\alpha\in(0,2)$, let $\Delta^{\frac{\alpha}{2}}$ be the fractional Laplacian in $\mR^d$ defined by
$$
\Delta^{\frac{\alpha}{2}}f(x)=
\lim_{\varepsilon \downarrow 0}\int_{|y|\geq\varepsilon}\frac{f(x+y)-f(x)}{|y|^{d+\alpha}}\dif y.
$$
It is well-known that the heat kernel $\rho^{(\alpha)}(t,x)$
of $\Delta^{\frac{\alpha}{2}}$ has the following estimate (e.g. see \cite{Ch-So, Ch-Ku}):
\begin{align}
\rho^{(\alpha)}(t,x)\asymp \frac{t}{(|x|\vee t^{\frac{1}{\alpha}})^{d+\alpha}},\label{Heat}
\end{align}
where $\asymp$ means that both sides are comparable up to some positive constants.

In \cite{Bo-Ja0}, Bogdan and Jakubowski studied the following perturbation of $\Delta^{\frac{\alpha}{2}}$
by gradient operator:
$$
\sL^{(\alpha)}_b(x):=\Delta^{\frac{\alpha}{2}}+b(x)\cdot\nabla,\ \ \alpha\in(1,2),
$$
where $b$ belongs to Kato's class $\sK^{\alpha-1}_d$, i.e., for $\gamma>0$,
$$
b\in\sK^{\gamma}_d:=\left\{f\in L^1_{loc}(\mR^d): \lim_{\eps\downarrow 0}\sup_{x\in\mR^d}
\int_{|x-y|\leq\eps}\frac{|f(y)|}{|x-y|^{d-\gamma}}\dif y=0\right\}.
$$
Notice that by H\"older's inequality, $L^p(\mR^d)\subset\sK^{\gamma}_d$ provided $p>\frac{d}{\gamma}$.
The sharp two-sided  heat kernel estimates of $\sL^{(\alpha)}_b$ like (\ref{Heat}) were obtained therein.
The reason of limiting $\alpha\in(1,2)$ lies in that the heat kernel $p^{(\alpha)}_1(t,x)=\rho^{(\alpha)}(t,x+t)$
of $\sL^{(\alpha)}_1$ is not comparable with $\rho^{(\alpha)}(t,x)$ for $\alpha\in(0,1)$ (see \cite{Bo-Ja0}).
In \cite{Ja-Sz2}, Jakubowski and Szczypkowski considered the time-dependent perturbation of $\Delta^{\frac{\alpha}{2}}$.
In \cite{Ja}, Jakubowski established the global time estimate of heat kernel of $\Delta^{\frac{\alpha}{2}}$ with small singular drift.
In \cite{Ch-Ki-So}, Chen, Kim and Song obtained sharp two-sided estimates for the Dirichlet heat kernel
of $\sL^{(\alpha)}_b$. Moreover, the Dirichlet heat kernel estimates for nonlocal operators under
Feynman-Kac or Schr\"odinger type perturbations were also considered in \cite{Ch-Ki-So1}.
Recently, in \cite{Wa-Zh}, Wang and the second named author extended Bogdan and
Jakubowski's results to the more general subordinated
stable operator over Riemannian manifold and obtained sharp two-sided estimates as well as the gradient estimate.

However, in the critical case of $\alpha=1$, the heat kernel estimate of $\sL^{(1)}_b$ is left open. It is noticed that
the critical case has particular interest in physics and mathematics (see \cite{Ca-Va,Ki-Na-Vo,Ki-Na-Sh,Si3,Si} and references therein).
We first recall some related results. In \cite{Ma-Mi}, Maekawa and Miura obtained the upper bounds estimates
for the fundamental solutions of general nonlocal diffusions with divergence free drift.
Their proofs are based upon the classical Davies' method. In \cite{Si3} and \cite{Si}, Silvestre established the H\"older regularity to the critical parabolic
operator $\sL^{(1)}_b(x)$ with bounded measurable $b$. In \cite{Pr},
Priola proved the pathwise uniqueness of SDEs with H\"older's drifts and driven by Cauchy processes.
In \cite{Zh1}, the well-posedness of multidimensional critical Burgers' equation was obtained (see \cite{Ki-Na-Sh}
for the study of one dimensional critical Burgers' equations).

In this paper we consider the following critical fractional diffusion operator
$$
\sL(t,x):=\sL_{a,b,c}(t,x):=a(t,x)\Delta^{\frac{1}{2}}+b(t,x)\cdot\nabla+c(t,x),
$$
where $a, c:[0,\infty)\times\mR^d\to\mR$ and $b:[0,\infty)\times\mR^d\to\mR^d$ are three measurable functions.
We shall prove the following result in the present work.
\bt\label{Main}
Assume that for some $a_0,a_1>0$,
$$
a_0\leq a(t,x)\leq a_1,
$$
and for some $\beta\in(0,1)$,
$$
a,b\in \mH^\beta,\ \ c\in\mK^1_d,
$$
where $\mH^\beta$ (resp. $\mK^1_d$) is the H\"older space (resp. the generalized Kato's class) defined in Definition \ref{Def1}.
Then there exists a continuous function $p(t,x;s,y)$ such that:
\begin{enumerate}[(i)]

\item (C-K equation) For all $0\leq s<r<t$ and $x,y\in\mR^d$, the following Chapman-Kolmogorov's equation holds:
\begin{align}\label{eq21}
\int_{\mR^d}p(t,x;r,z)p(r,z;s,y)\dif z=p(t,x;s,y).
\end{align}

\item (Generator) For any bounded continuous function $f$, we have
\begin{align}\label{eq22}
\lim_{t\downarrow s}P_{t,s}f(x):=\lim_{t\downarrow s}\int_{\mR^d}p(t,x;s,y)f(y)\dif y=f(x),
\end{align}
and if $a,b,c\in C\big([0,\infty);L_{\mathrm{loc}}^{1}(\mR^d)\big)$, then for all $f,g\in C^\infty_0(\mR^d)$,
\begin{align}\label{eq23}
\lim_{t\downarrow s}\frac{1}{t-s}\int_{\mR^d}g(x)\Big(P_{t,s}f(x)-f(x)\Big)\dif x
=\int_{\mR^d}g(x)\sL(s,x)f(x)\dif x,\ \ s\geq 0.
\end{align}

\item (Two-sided estimates) For any $T>0$, there exist constants
$\kappa_1,\kappa_2>0$ such that for all $0\leq s<t\leq T$ and $x,y\in\mR^d$,
\begin{align}\label{eq24}
|p(t,x;s,y)|\leq \kappa_1(t-s)(|x-y|+(t-s))^{-d-1},
\end{align}
and in the case of $a(t,x)=a(t)$ independent of $x$,
\begin{align}\label{eq25}
p(t,x;s,y)\geq \kappa_2(t-s)(|x-y|+(t-s))^{-d-1}.
\end{align}
\item (H\"older's estimate) Assume that $c\in\mK^{1-\gamma}_d$ for some $\gamma\in(0,1)$.
Then for any $T>0$, there exists a constant $\kappa_3>0$ such that
 for all $0\leq s<t\leq T$ and $x,x',y\in\mR^d$,
\begin{align}
|p(t,x;s,y)-p(t,x';s,y)|&\leq \kappa_3(|x-x'|^{\gamma}\wedge 1)|t-s|^{1-\gamma}\no\\
&\times\Big\{(|x-y|+(t-s))^{-d-1}
+(|x'-y|+(t-s))^{-d-1}\Big\}.\label{eq34}
\end{align}
\item(Gradient estimate) If we further assume that $c\in\mH^{\gamma}$ for some $\gamma\in(0,1)$,
then for any $T>0$, there exists a constant $\kappa_4>0$ such that  for all $0\leq s<t\leq T$ and $x,y\in\mR^d$,
\begin{align}\label{eq26}
|\nabla_x p(t,x;s,y)|\leq \kappa_4(|x-y|+(t-s))^{-d-1}.
\end{align}
\end{enumerate}
\et
In order to prove this theorem, we shall use Levi's method of freezing coefficients and Duhamel's formula.
Compared with the classical case of second order parabolic equations, the main difficulty of proving this theorem
lies in the heavy tail property of Poisson's kernel and the nonlocal property of $\Delta^{\frac{1}{2}}$.
We mention that in the case of second order parabolic equation, the following property of Gaussian heat kernel plays a key role
in the construction of Levi's argument (cf. \cite{Fr, La-So-Ur}):
for $\beta\in(0,1)$ and some $C>0$,
$$
t^{-1}|x|^\beta \mathrm{e}^{-\frac{|x|^2}{t}}\leq t^{\frac{\beta}{2}-1}\mathrm{e}^{-\frac{|x|^2}{Ct}},\ \ t>0, x\in\mR^d.
$$
This means that the spatial H\"older regularity can compensate the time singularity. However, such type estimate
does not hold for Poisson's kernel in view of the heavy tail property. A suitable
substitution is an analogue of the so called 3P-inequality (see Lemma \ref{Le2} below).

This paper is organized as follows: In Section 2, we prepare some lemmas for later use.
In Section 3, by using Levi's method of constructing the fundamental solutions,
we first construct the heat kernel of $\sL_{a,b}=\sL_{a,b,0}$. In Section 4, we prove Theorem \ref{Main}
by using Duhamel's formula.

We conclude this section by introducing the following conventions: The letter $C$ with or without subscripts
will denote a positive constant, whose value is not important and may change in different places. We write $f(x)\preceq g(x)$ to mean that there exists a constant
$C_0>0$ such that $f(x)\leq C_0 g(x)$; and $f(x)\asymp g(x)$ to mean
that there exist $C_1,C_2>0$ such that $C_1 g(x)\leq f(x)\leq C_2 g(x)$.

\section{Preliminaries}

For $\gamma,\beta\in\mR$, we introduce the following function on $\mR_+\times\mR^d$ for later use:
\begin{align}
\varrho^\beta_\gamma(t,x):=t^\gamma\{|x|^\beta\wedge1\}(|x|^2+t^2)^{-\frac{d+1}{2}}
\asymp t^\gamma\{|x|^\beta\wedge1\}(|x|+t)^{-d-1}.\label{Def2}
\end{align}
By simple calculations, there exists a constant $C_d>0$ such that for all $\beta\in[0,\frac{1}{2}]$ and $\gamma\in\mR$,
\begin{align}
\int_{\mR^d}\varrho^\beta_\gamma(t,x)\dif x\leq C_d t^{\gamma+\beta-1}.\label{ES4}
\end{align}
Indeed,  we have
\begin{align*}
\int_{\mR^d}\frac{|x|^\beta}{(|x|+t)^{d+1}}\dif x
&\preceq\int^\infty_0\frac{r^{\beta+d-1}}{(r+t)^{d+1}}\dif r
=\left(\int^t_0+\int^\infty_t\right)\frac{r^{\beta+d-1}}{(r+t)^{d+1}}\dif r\\
&\leq \int^t_0\frac{r^{d+\beta-1}}{t^{d+1}}\dif r+\int^\infty_tr^{\beta-2}\dif r
=\frac{t^{\beta-1}}{d+\beta}+\frac{t^{\beta-1}}{1-\beta},
\end{align*}
which in turn implies (\ref{ES4}). Notice that the following $3$P-inequality holds (cf. \cite[Lemma 2.1]{Bo-Ja0}):
\begin{align}
\varrho^0_1(t,x)\varrho^0_1(s,y)\preceq\big(\varrho^0_1(t,x)+\varrho^0_1(s,y)\big)\varrho^0_1(t+s,x+y).\label{eq6}
\end{align}

We introduce the following classes of functions used  in this paper.
\bd\label{Def1}
(H\"older's space) For $\beta\in(0,1]$, define
$$
\mH^\beta:=\left\{f\in\sB(\mR\times\mR^d): \|f\|_{\mH^\beta}:=\sup_{t\in\mR}\sup_{x\in\mR^d}|f(t,x)|
+\sup_{t\in\mR}\sup_{x\not= y\in\mR^d}\frac{|f(t,x)-f(t,y)|}{|x-y|^\beta}<\infty\right\}.
$$
(Generalized Kato's class) For $\gamma>0$, define
$$
\mK^\gamma_d:=\Big\{f\in L^1_{\mathrm{loc}}(\mR\times\mR^d): \lim_{\eps\downarrow 0}K^\gamma(\eps)=0\Big\},
$$
where
$$
K^\gamma(\eps):=\sup_{(t,x)\in[0,\infty)\times\mR^d} \int^\eps_0\!\!\!\int_{\mR^d}
\varrho^0_\gamma(s,y)|f(t\pm s,x-y)| \dif y\dif s,\ \ \eps>0.
$$
\ed
A function $f(t,x)$ on $[0,\infty)\times\mR^d$ will be automatically extended to $\mR\times\mR^d$ by letting $f(t,\cdot)=0$ for $t<0$.
The following proposition gives a characterization for $\mK^\gamma_d$ (see \cite{Ai-Si, Zh0, Wa-Zh} for more discussions).
\bp
For $\gamma>0$ and $p,q\in[1,\infty]$ with $\frac{d}{p}+\frac{1}{q}<\gamma$, we have
$$
L^q(\mR;L^p(\mR^d))\subset\mK^\gamma_d,
$$
and for $\gamma\in(0,d)$,
$$
\sK^\gamma_d\subset\mK^\gamma_d.
$$
\ep
\begin{proof}
Noticing that
\begin{align*}
\int^\eps_0\!\!\!\int_{\mR^d}\varrho^0_\gamma(s,y)|f(t\pm s,x-y)| \dif y\dif s=
\int^\eps_0\!\!s^{\gamma-1}\!\!\!\int_{\mR^d}\varrho^0_1(s,y)|f(t\pm s,x-y)| \dif y\dif s,
\end{align*}
by H\"older's inequality, for the first inclusion, it is enough to prove
\begin{align}
\lim_{\eps\downarrow 0}I(\eps)=0,\label{eq0}
\end{align}
where
$$
I(\eps):=\int^\eps_0\!\Bigg(\!\int_{\mR^d}{\varrho^0_{1}(s,y)}^{p^{*}}\dif y
\Bigg)^{\frac{q^*}{p^*}}\!\!s^{(\gamma-1) q^*}\dif s
$$
with $q^*=\frac{q}{q-1}$ and $p^*=\frac{p}{p-1}$. As in the proof of (\ref{ES4}), we have
$$
\int_{\mR^d}{\varrho^0_{1}(s,y)}^{p^{*}}\dif y\preceq s^{d-dp^*},
$$
and so,
$$
I(\eps)\preceq\int^\eps_0 s^{\frac{dq^*}{p^*}-dq^*+(\gamma-1)q^*}\dif s\preceq\eps^{1+\frac{dq^*}{p^*}-dq^*+(\gamma-1) q^*},
$$
since $\frac{dq^*}{p^*}-dq^*+(\gamma-1) q^*>-1$ by $\frac{d}{p}+\frac{1}{q}<\gamma$, thus (\ref{eq0}) holds.

Next we prove the second inclusion. Assume $f\in\sK^\gamma_d$. By definition, we have
$$
\sup_{x\in\mR^d}\int^\eps_0\!\!\!\int_{\mR^d}\varrho^0_\gamma(s,y)|f(x-y)|\dif y\dif s\leq I_1(\eps)+I_2(\eps),\
$$
where
\begin{align*}
I_1(\eps)&:=\sup_{x\in\mR^d} \int^\eps_0\!\!\!\int_{|y|\leq\eps} \frac{s^\gamma|f(x-y)|}{(|y|+s)^{d+1}}\dif y\dif s,\\
I_2(\eps)&:=\sup_{x\in\mR^d} \int^\eps_0\!\!\!\int_{|y|>\eps} \frac{s^\gamma|f(x-y)|}{(|y|+s)^{d+1}}\dif y\dif s.
\end{align*}
For $I_1(\eps)$, in view of $\gamma<d$, we have
\begin{align*}
I_1(\eps)&\leq\sup_{x\in\mR^d} \int_{|y|\leq\eps} |f(x-y)|
\left(\int^\eps_{|y|}s^{\gamma-d-1}\dif s+|y|^{-d-1}\int^{|y|}_0s^\gamma\dif s\right)\dif y\\
&\leq\sup_{x\in\mR^d} \int_{|y|\leq\eps}|f(x-y)|\left(\frac{|y|^{-d+\gamma}}{d-\gamma}
+\frac{|y|^{-d+\gamma}}{\gamma+1}\right)\dif y\to 0,\ \ \eps\downarrow 0.
\end{align*}
For $I_2(\eps)$, we have
\begin{align*}
I_2(\eps)&\leq\sup_{x\in\mR^d} \int_{|y|>\eps} \frac{|f(x-y)|}{|y|^{d+1}}\dif y\int^\eps_0s^\gamma\dif s
=\frac{1}{\gamma+1}\sup_{x\in\mR^d} \int_{|y|>\eps} \frac{\eps^{\gamma+1}|f(x-y)|}{|y|^{d+1}}\dif y,
\end{align*}
which converges to zero by \cite[Lemma 11]{Bo-Ja0} as $\eps\downarrow 0$.
\end{proof}
Set for $s<t$ and $x,y\in\mR^d$,
$$
\varrho^\beta_\gamma(t,x;s,y):=\varrho^\beta_\gamma(t-s,x-y).
$$
The following lemma is an analogue of $3$P-inequality, which will play a crucial role in the sequel.
\bl\label{Le2}
For all $\beta_1,\beta_2\in[0,\frac{1}{4}]$ and $\gamma_1,\gamma_2\in\mR$, we have
\begin{align}
\int_{\mR^d}\varrho^{\beta_1}_{\gamma_1}(t,x;r,z)\varrho^{\beta_2}_{\gamma_2}(r,z;s,y)\dif z
&\leq C_d\Big\{(t-r)^{\gamma_1+\beta_1+\beta_2-1}(r-s)^{\gamma_2}\varrho^0_0(t,x;s,y)\no\\
&\quad+(t-r)^{\gamma_1+\beta_1-1}(r-s)^{\gamma_2}\varrho^{\beta_2}_0(t,x;s,y)\no\\
&\quad+(t-r)^{\gamma_1}(r-s)^{\gamma_2+\beta_1+\beta_2-1}\varrho^0_0(t,x;s,y)\no\\
&\quad+(t-r)^{\gamma_1}(r-s)^{\gamma_2+\beta_2-1}\varrho^{\beta_1}_0(t,x;s,y)\Big\},\label{EU7}
\end{align}
where $C_d$ only depends on $d$; and if $\gamma_1>-\beta_1$ and $\gamma_2>-\beta_2$, then
\begin{align}
\int^t_s\!\!\!\int_{\mR^d}\varrho^{\beta_1}_{\gamma_1}(t,x;r,z)\varrho^{\beta_2}_{\gamma_2}(r,z;s,y)\dif z\dif r
&\leq C_d\Big\{\varrho^0_{\gamma_1+\gamma_2+\beta_1+\beta_2}(t,x;s,y)\cB(\gamma_1+\beta_1+\beta_2,1+\gamma_2)\no\\
&\quad+\varrho^{\beta_2}_{\gamma_1+\gamma_2+\beta_1}(t,x;s,y)\cB(\gamma_1+\beta_1,1+\gamma_2)\no\\
&\quad+\varrho^0_{\gamma_1+\gamma_2+\beta_1+\beta_2}(t,x;s,y)\cB(\gamma_2+\beta_1+\beta_2,1+\gamma_1))\no\\
&\quad+\varrho^{\beta_1}_{\gamma_1+\gamma_2+\beta_2}(t,x;s,y)\cB(\gamma_2+\beta_2,1+\gamma_1\Big\},\label{eq30}
\end{align}
where $\cB(\gamma,\beta)$ is the usual Beta function defined by
$$
\cB(\gamma,\beta):=\int^1_0(1-s)^{\gamma-1}s^{\beta-1}\dif s,\ \
\gamma,\beta>0.
$$
Moreover, there exist $p>1$ and a constant $C>0$ such that for all $0\leq s<t\leq 1$ and $x\not=y\in\mR^d$,
\begin{align}
\int^t_s\left(\int_{\mR^d}\varrho^{\beta_1}_{\gamma_1}(t,x;r,z)\varrho^{\beta_2}_{\gamma_2}
(r,z;s,y)\dif z\right)^p\dif r\leq \frac{C}{|x-y|^{(d+1)p}}.\label{eq31}
\end{align}
\el
\begin{proof}
First of all, in view of
$$
(|x-y|^2+|t-s|^2)^{\frac{d+1}{2}}\leq
2^d\Big\{(|x-z|^2+|t-r|^2)^{\frac{d+1}{2}}+(|z-y|^2+|r-s|^2)^{\frac{d+1}{2}}\Big\},
$$
we have
\begin{align}
\varrho^0_0(t,x;r,z)\varrho^0_0(r,z;s,y)\leq
2^d\Big(\varrho^0_0(t,x;r,z)+\varrho^0_0(r,z;s,y)\Big)\varrho^0_0(t,x;s,y).\label{ET2}
\end{align}
Noticing that by $(a+b)^\beta\leq a^\beta+b^\beta$ for $\beta\in(0,1)$,
\begin{align*}
(|x-z|^{\beta_1}\wedge 1)(|z-y|^{\beta_2}\wedge 1)
&\leq(|x-z|^{\beta_1}\wedge 1)((|x-z|^{\beta_2}+|x-y|^{\beta_2})\wedge 1)\\
&\leq|x-z|^{\beta_1+\beta_2}\wedge 1+(|x-z|^{\beta_1}\wedge 1)(|x-y|^{\beta_2}\wedge 1),\\
(|x-z|^{\beta_1}\wedge 1)(|z-y|^{\beta_2}\wedge 1)
&\leq((|z-y|^{\beta_1}+|x-y|^{\beta_1})\wedge 1)(|z-y|^{\beta_2}\wedge 1)\\
&\leq|z-y|^{\beta_1+\beta_2}\wedge 1+(|z-y|^{\beta_2}\wedge 1)(|x-y|^{\beta_1}\wedge1),
\end{align*}
we have
\begin{align*}
&\varrho^{\beta_1}_{\gamma_1}(t,x;r,z)\varrho^{\beta_2}_{\gamma_2}(r,z;s,y)\no\\
&\qquad=|t-r|^{\gamma_1}|r-s|^{\gamma_2}(|x-z|^{\beta_1}\wedge 1)(|z-y|^{\beta_2}\wedge 1)
\varrho^0_0(t,x;r,z)\varrho^0_0(r,z;s,y)\no\\
&\qquad\preceq|t-r|^{\gamma_1}|r-s|^{\gamma_2}\Big((|x-z|^{\beta_1}\wedge 1)(|x-y|^{\beta_2}\wedge 1)
+|x-z|^{\beta_1+\beta_2}\wedge 1\Big)\no\\
&\qquad\qquad\qquad\qquad\quad\times\varrho^0_0(t,x;r,z)\varrho^0_0(t,x;s,y)\no\\
&\qquad+|t-r|^{\gamma_1}|r-s|^{\gamma_2}\Big((|z-y|^{\beta_2}\wedge 1)(|x-y|^{\beta_1}\wedge 1)
+|z-y|^{\beta_1+\beta_2}\wedge 1\Big)\no\\
&\qquad\qquad\qquad\qquad\quad\times\varrho^0_0(r,z;s,y)\varrho^0_0(t,x;s,y)\no\\
&\qquad\preceq |r-s|^{\gamma_2}(\varrho^{\beta_1+\beta_2}_{\gamma_1}(t,x;r,z)\varrho^0_0(t,x;s,y)
+\varrho^{\beta_1}_{\gamma_1}(t,x;r,z)\varrho^{\beta_2}_0(t,x;s,y))\no\\
&\qquad+|t-r|^{\gamma_1}(\varrho^{\beta_1+\beta_2}_{\gamma_2}(r,z;s,y)\varrho^0_0(t,x;s,y)
+\varrho^{\beta_2}_{\gamma_2}(r,z;s,y)\varrho^{\beta_1}_0(t,x;s,y)).
\end{align*}
Estimate (\ref{EU7}) follows by (\ref{ES4}), and estimate (\ref{eq30}) follows by observing that for $\gamma,\beta>0$,
\begin{align}
\int^t_s(t-r)^{\gamma-1}(r-s)^{\beta-1}\dif r=(t-s)^{\gamma+\beta-1}\cB(\gamma,\beta).\label{eq32}
\end{align}
As for estimate (\ref{eq31}), it follows by (\ref{EU7}) and (\ref{eq32}).
\end{proof}

Let $\rho(t,x)$ be the heat kernel of the Cauchy operator $\Delta^{\frac{1}{2}}$, i.e.,
\begin{align}
\p_t\rho(t,x)=\Delta^{\frac{1}{2}}\rho(t,x).\label{ES1}
\end{align}
It is well-known that
$$
\rho(t,x)=\pi^{-\frac{d+1}{2}}\Gamma(\tfrac{d+1}{2})(|x|^2+t^2)^{-\frac{d+1}{2}}t
=\pi^{-\frac{d+1}{2}}\Gamma(\tfrac{d+1}{2})\varrho^0_1(t,x),
$$
which is also called Poisson kernel (cf. \cite{St}). By elementary calculations, one has
\begin{align}
|\nabla_x\rho(t,x)|\preceq t(|x|+t)^{-d-2},\ & \ |\p_t\rho(t,x)|\preceq (|x|+t)^{-d-1},\label{EU1}\\
|\nabla^2_x\rho(t,x)|+|\nabla_x\p_t\rho(t,x)|&\preceq (|x|+t)^{-d-2},\label{EU3}
\end{align}
and
\begin{align}
|\nabla^3_x\rho(t,x)|+|\nabla^2_x\p_t\rho(t,x)|\preceq (|x|+t)^{-d-3}.\label{EU33}
\end{align}

Let $a: [0,\infty)\times\mR^d\to(0,\infty)$ and $b:[0,\infty)\times\mR^d\to\mR^d$ be two bounded measurable functions.
We define
$$
p_0(t,x;s,y):=\rho\left(\int^t_s a(r,y)\dif r,x-y+\int^t_sb(r,y)\dif r\right),
$$
and
\begin{align}
\sL^x_{a,b}(t,y):=a(t,y)\Delta^{\frac{1}{2}}_x+b(t,y)\cdot\nabla_x.
\end{align}
By (\ref{ES1}) and the Lebesgue differential theorem, we have for all $x,y\in\mR^d$ and almost all $t>s$,
\begin{align}
\p_tp_0(t,x;s,y)=\sL^x_{a,b}(t,y)p_0(t,\cdot;s,y)(x).\label{ES2}
\end{align}

We prepare the following important estimates for later use.
\bl\label{Le1} Suppose that for some $a_0,a_1,b_1>0$,
\begin{align}
a_0\leq a(r,y)\leq a_1,\ \ |b(r,y)|\leq b_1.\label{Con}
\end{align}
Then we have
\begin{align}
p_0(t,x;s,y)\asymp\varrho^0_1(t,x;s,y),\label{E1}
\end{align}
and
\begin{align}
|\Delta^{\frac{1}{2}}_xp_0(t,x;s,y)|&\preceq \big(|x-y|+|t-s|\big)^{-d-1},\label{EU5}\\
|\nabla_xp_0(t,x;s,y)|&\preceq |t-s|\big(|x-y|+|t-s|\big)^{-d-2},\label{EU4}\\
|\p_tp_0(t,x;s,y)|&\preceq  \big(|x-y|+|t-s|\big)^{-d-1},\label{EU404}\\
|\nabla_x\Delta^{\frac{1}{2}}_xp_0(t,x;s,y)|&\preceq \big(|x-y|+|t-s|\big)^{-d-2},\label{EU6}\\
|\nabla^2_xp_0(t,x;s,y)|&\preceq\big(|x-y|+|t-s|\big)^{-d-2}.\label{EU44}
\end{align}
Moreover, if we further assume that $a,b\in\mH^\beta$ for some $\beta\in(0,1)$, then
\begin{align}
\left|\int_{\mR^d}\nabla_xp_0(t,x;s,y)\dif y\right|&\preceq(t-s)^{\beta-1},\label{EW7}\\
\left|\int_{\mR^d}\Delta^{\frac{1}{2}}_xp_0(t,x;s,y)\dif y\right|&\preceq(t-s)^{\beta-1},\label{EW8}\\
\left|\int_{\mR^d}\p_tp_0(t,x;s,y)\dif y\right|&\preceq(t-s)^{\beta-1},\label{EW6}\\
\lim_{t\downarrow s}\int_{\mR^d}p_0(t,x;s,y)\dif y&=1,\label{EW66}
\end{align}
and for all $w\in\mR^d$ and $\gamma\in[0,\beta]$,
\begin{align}
\left|\int_{\mR^d}(\nabla_xp_0(t,x+w;s,y)-\nabla_xp_0(t,x;s,y))\dif y\right|\preceq
|w|^{\gamma}(t-s)^{\beta-\gamma-1}.\label{EQ1}
\end{align}
\el
\begin{proof}
For the simplicity of notation, we write
$$
F^t_s(y):=\int^t_s a(r,y)\dif r,\ \ G^t_s(y):=\int^t_s b(r,y)\dif r.
$$

(1) By (\ref{Con}), we have
\begin{align}
F^t_s(y)\asymp t-s,\label{EQ33}
\end{align}
and for any $|w|\leq|t-s|$,
\begin{align}
\left|x+w-y+G^t_s(y)\dif r\right|+|t-s|\asymp |x-y|+|t-s|.\label{EQ3}
\end{align}
Estimate (\ref{E1}) follows by definition. For (\ref{EU5}), by (\ref{ES1})
we have
\begin{align*}
\Delta^{\frac{1}{2}}_xp_0(t,x;s,y)=(\Delta^{\frac{1}{2}}_x\rho)
\left(F^t_s(y), x-y+G^t_s(y)\right)=(\p_t\rho)\left(F^t_s(y), x-y+G^t_s(y)\right).
\end{align*}
Estimate (\ref{EU5}) follows by (\ref{EU1}). Similarly,  (\ref{EU4})-(\ref{EU44})
follow by (\ref{EU1}), (\ref{EU3})  and  (\ref{ES2}).

(2) Define
$$
\xi(t,x;s,y;z):=\rho\left(\int^t_s a(r,z)\dif r,x-y+\int^t_sb(r,z)\dif r\right)=\rho\Big(F^t_s(z),x-y+G^t_s(z)\Big).
$$
Clearly, for any $s<t$ and $x,z\in\mR^d$,
$$
\int_{\mR^d}\xi(t,x;s,y;z)\dif y=\int_{\mR^d}\rho\left(F^t_s(z),y\right)\dif y=1
$$
and
\begin{align*}
\int_{\mR^d}\nabla_x\xi(t,x;s,y;z)\dif y=0,\ \ \
\int_{\mR^d}\Delta^{\frac{1}{2}}_x\xi(t,x;s,y;z)\dif y=0.
\end{align*}
Thus, for proving (\ref{EW7}), it suffices to prove that
\begin{align}
\left|\int_{\mR^d}\Big(\nabla_x p_0(t,x;s,y)-\nabla_x\xi(t,x;s,y;z)\Big)\dif y\right|_{z=x}\preceq(t-s)^{\beta-1}.\label{EP3}
\end{align}
By $a,b\in\mH^\beta$ and definitions of $p_0$ and $\xi$, one has
\begin{align}
&|\nabla_x p_0(t,x;s,y)-\nabla_x\xi(t,x;s,y;z)|_{z=x}|\no\\
&\quad=|(\nabla_x\rho)\left(F^t_s(y),x-y+G^t_s(y)\right)
-(\nabla_x\rho)\left(F^t_s(x),x-y+G^t_s(x)\right)|\no\\
&\quad\leq\|a\|_{\mH^\beta}(|x-y|^\beta\wedge 1)|t-s|\int^1_0|\nabla_x\p_t\rho|\left(\theta F^t_s(y)
+(1-\theta)F^t_s(x),x-y+G^t_s(y)\right)\dif\theta\no\\
&\quad+\|b\|_{\mH^\beta}(|x-y|^\beta\wedge 1)|t-s|\int^1_0|\nabla_x^2\rho|\left(F^t_s(x),x-y+\theta G^t_s(y)
+(1-\theta)G^t_s(x)\right)\dif\theta\no\\
&\quad\stackrel{(\ref{EU3}),(\ref{EQ33}),(\ref{EQ3})}\preceq \frac{(|x-y|^\beta\wedge 1)|t-s|}
{(|x-y|+|t-s|)^{d+2}}\leq\frac{(|x-y|^\beta\wedge 1)}{(|x-y|+|t-s|)^{d+1}},\label{EQ2}
\end{align}
which then gives estimate (\ref{EP3}) by (\ref{ES4}).

Similarly, we can prove
$$
\left|\int_{\mR^d}\Big(\Delta^{\frac{1}{2}}_x p_0(t,x;s,y)
-\Delta^{\frac{1}{2}}_x\xi(t,x;s,y;z)\Big)\dif y\right|_{z=x}\preceq(t-s)^{\beta-1},
$$
and
$$
\left|\int_{\mR^d}\Big(p_0(t,x;s,y)-\xi(t,x;s,y;z)\Big)\dif y\right|_{z=x}\preceq(t-s)^{\beta}.
$$
Thus, (\ref{EW8}) and (\ref{EW66}) follow.

(3) Next, we prove (\ref{EW6}). By (\ref{ES2}), (\ref{EU5}), (\ref{EU4}), (\ref{EW7}) and (\ref{EW8}), we have
\begin{align*}
\left|\int_{\mR^d}\p_tp_0(t,x;s,y)\dif y\right|
&=\left|\int_{\mR^d}\Big(a(t,y)\Delta^{\frac{1}{2}}_xp_0(t,x;s,y)+b(t,y)\cdot\nabla_xp_0(t,x;s,y)\Big)\dif y\right|\\
&\leq|a(t,x)|\left|\int_{\mR^d}\Delta^{\frac{1}{2}}_xp_0(t,x;s,y)\dif y\right|
+|b(t,x)|\left|\int_{\mR^d}\nabla_xp_0(t,x;s,y)\dif y\right|\\
&\quad+\int_{\mR^d}|a(t,y)-a(t,x)|\cdot|\Delta^{\frac{1}{2}}_xp_0(t,x;s,y)|\dif y\\
&\quad+\int_{\mR^d}|b(t,y)-b(t,x)|\cdot|\nabla_xp_0(t,x;s,y)|\dif y\\
&\preceq(t-s)^{\beta-1}+\int_{\mR^d}\varrho^\beta_0(t,x;s,y)\dif y\stackrel{(\ref{ES4})}{\preceq}(t-s)^{\beta-1}.
\end{align*}

(4) Lastly, we prove (\ref{EQ1}). 
If $|w|\leq|t-s|$, then
\begin{align*}
&|\nabla_x p_0(t,x+w;s,y)-\nabla_x\xi(t,x+w;s,y;z)|_{z=x}
-(\nabla_x p_0(t,x;s,y)-\nabla_x\xi(t,x;s,y;z)|_{z=x})|\\
&=\Bigg|w\cdot\int^1_0\Big[(\nabla^2_x\rho)\left(F^t_s(y),x+\theta w-y+G^t_s(y)\right)
-(\nabla^2_x\rho)\left(F^t_s(x),x+\theta w-y+F^t_s(x)\dif r\right)\Big]\dif \theta\Bigg|\\
&\preceq|w|\frac{(t-s)(|x-y|^\beta\wedge 1)}{(|x-y|+(t-s))^{d+3}}
\preceq|w|^{\gamma}(t-s)^{\beta-\gamma}\varrho^0_0(t,x;s;y),
\end{align*}
where we have used the same argument as in proving (\ref{EQ2}). Integrating both sides with respect to $y$ and using (\ref{ES4}),
we obtain (\ref{EQ1}) for $|w|\leq|t-s|$. If $|w|>|t-s|$, it follows by (\ref{EW7}).
\end{proof}
\br
By (\ref{EU4}), we also have for any $\gamma\in(0,1]$,
\begin{align}
|\nabla_xp_0(t,x;s,y)|\preceq |t-s|^\gamma|x-y|^{-\gamma}\varrho^0_0(t,x;s,y).\label{EW5}
\end{align}
This estimate is important for the lower bound estimate of the heat kernel.
\er

\section{Heat kernel of $\sL_{a,b}:=a\Delta^{\frac{1}{2}}+b\cdot\nabla$}

Let $\sL_{a,b}:=\sL^x_{a,b}(t,x)=a(t,x)\Delta^{\frac{1}{2}}_x+b(t,x)\cdot\nabla_x$.
Now we want to seek the heat kernel of $\sL_{a,b}$ with the following form:
\begin{align}
p_{a,b}(t,x;s,y)=p_0(t,x;s,y)+\int^t_s\!\!\!\int_{\mR^d}p_0(t,x;r,z)q(r,z;s,y)\dif z\dif r. \label{eq2}
\end{align}
The classical Levi's continuity argument (see \cite{La-So-Ur,Fr}) suggests that $q(t,x;s,y)$ must satisfy the following integral equation:
\begin{align}
q(t,x;s,y)=q_0(t,x;s,y)+\int^t_s\!\!\!\int_{\mR^d}q_0(t,x;r,z)q(r,z;s,y)\dif z\dif r,\label{Eq}
\end{align}
where
\begin{align}
q_0(t,x;s,y):=\big(a(t,x)-a(t,y)\big)\Delta^{\frac{1}{2}}_xp_0(t,x;s,y)+\big(b(t,x)-b(t,y)\big)\cdot\nabla_xp_0(t,x;s,y).
\label{eq1}
\end{align}
In the remainder of this paper, we shall work on the time interval $[0,1]$, and always assume
$$
0\leq s<t\leq 1,\ \ x\not=y\in\mR^d,
$$
and for some $\beta\in(0,1)$,
\begin{align}
a,b\in \mH^\beta.\label{eq40}
\end{align}
Our first task is thus to solve the integral equation (\ref{Eq}).

Let us now recursively define
\begin{align}
q_n(t,x;s,y):=\int^t_s\!\!\!\int_{\mR^d}q_0(t,x;r,z)q_{n-1}(r,z;s,y)\dif z\dif r,\ \ n\in\mN.\label{ES6}
\end{align}

\bl\label{Le3}
For $\beta\in(0,\frac{1}{4}]$, there exists a constant $C_d>0$ such that for all $n\in\mN$,
\begin{align}
|q_n(t,x;s,y)|\leq\frac{(C_d\Gamma(\beta))^{n+1}}{\Gamma((n+1)\beta)}\left(\varrho^0_{(n+1)\beta}(t,x;s,y)
+\varrho^\beta_{n\beta}(t,x;s,y)\right).\label{EW2}
\end{align}
Moreover, if $a(t,x)=a(t)$ is independent of $x$, then
\begin{align}
|q_n(t,x;s,y)|\leq\frac{(C_d\Gamma(\beta))^{n+1}}{\Gamma((n+1)\beta)}\varrho^0_{(n+1)\beta}(t,x;s,y).\label{EW3}
\end{align}
\el
\begin{proof}
First of all, by (\ref{eq40}) and Lemma \ref{Le1}, we have
$$
|q_0(t,x;s,y)|\leq C_d\varrho^\beta_0(t,x;s,y).
$$
Notice that
$$
\cB(\gamma,\beta) \mbox{ is symmetric and non-increasing with respect to each
variable $\gamma$ and $\beta$}.
$$
For $n=1$, by Lemma \ref{Le2}, we have
$$
|q_1|\leq C_d\cB(2\beta,1)\varrho^0_{2\beta}+C_d\cB(\beta,1)\varrho^\beta_\beta
\leq C_d\cB(\beta,\beta)\Big\{\varrho^0_{2\beta}+\varrho^\beta_\beta\Big\}.
$$
Suppose now that
$$
|q_n|\leq\gamma_n\Big\{\varrho^0_{(n+1)\beta}+\varrho^\beta_{n\beta}\Big\},
$$
where $\gamma_n>0$ will be determined below.
By Lemma \ref{Le2} we have
\begin{align*}
|q_{n+1}|&\leq C_d\gamma_n\Big\{\cB(\beta,1+(n+1)\beta)+\cB((n+2)\beta,1)+\cB(2\beta,1+n\beta)\Big\}\varrho^0_{(n+2)\beta}\\
&\quad+C_d\gamma_{n} \Big\{\cB((n+1)\beta,1)+\cB(\beta,1+n\beta)\Big\}\varrho^\beta_{(n+1)\beta}\\
&\leq C_d\gamma_{n}
\cB(\beta,(n+1)\beta)\Big\{\varrho^0_{(n+2)\beta}+\varrho^\beta_{(n+1)\beta}\Big\}
=:\gamma_{n+1}\Big\{\varrho^0_{(n+2)\beta}+\varrho^\beta_{(n+1)\beta}\Big\},
\end{align*}
where
$$
\gamma_{n+1}=C_d\gamma_{n} \cB(\beta,(n+1)\beta).
$$
Hence, by
$\cB(\gamma,\beta)=\frac{\Gamma(\gamma)\Gamma(\beta)}{\Gamma(\gamma+\beta)}$,
we obtain
\begin{align*}
\gamma_{n}= C_d^{n+1}\cB(\beta,\beta)\cB(\beta,2\beta)\cdots\cB(\beta,n\beta)=
\frac{(C_d\Gamma(\beta))^{n+1}}{\Gamma((n+1)\beta)},
\end{align*}
which gives (\ref{EW2}).

In the case of $a(t,x)=a(t)$, by (\ref{EW5}) we have
$$
|q_0(t,x;s,y)|\preceq \varrho^0_\beta(t,x;s,y).
$$
Repeating the above proof, we obtain (\ref{EW3}).
\end{proof}

We also need the following H\"older continuity of $q_n$ with respect to $x$.
\bl\label{Le4} For all $n\geq 0$ and $\gamma\in(0,\beta)$, we have
\begin{align*}
&|q_n(t,x;s,y)-q_n(t,x';s,y)|\preceq \frac{(C_d\Gamma(\beta))^{n+1}}{\Gamma(n\beta+\gamma)}\Big(|x-x'|^{\beta-\gamma}\wedge 1\Big)\\
&\quad\times\Bigg\{\Big(\varrho^0_{\gamma+n\beta}+\varrho^\beta_{\gamma+(n-1)\beta}\Big)(t,x;s,y)
+\Big(\varrho^0_{\gamma+n\beta}+\varrho^\beta_{\gamma+(n-1)\beta}\Big)(t,x';s,y)\Bigg\}.
\end{align*}
\el
\begin{proof}
Let us first prove the following estimate:
\begin{align}
&|q_0(t,x;s,y)-q_0(t,x';s,y)|\no\\
&\qquad\preceq(|x-x'|^{\beta-\gamma}\wedge 1)
\Big\{(\varrho^0_\gamma+\varrho^\beta_{\gamma-\beta})(t,x;s,y)
+(\varrho^0_\gamma+\varrho^\beta_{\gamma-\beta})(t,x';s,y)\Big\}.\label{ES5}
\end{align}
In the case of $|x-x'|>1$, we have
$$
|q_0(t,x;s,y)|\preceq(\varrho^0_\beta+\varrho^\beta_0)(t,x;s,y)
\leq(\varrho^0_\gamma+\varrho^\beta_{\gamma-\beta})(t,x;s,y)
$$
and
$$
|q_0(t,x';s,y)|\preceq(\varrho^0_\beta+\varrho^\beta_0)(t,x';s,y)
\leq(\varrho^0_\gamma+\varrho^\beta_{\gamma-\beta})(t,x';s,y).
$$
In the case of $1\geq |x-x'|>|t-s|$, by (\ref{EU5}) and (\ref{EU4}) we have
$$
|q_0(t,x;s,y)|\preceq\varrho^\beta_0(t,x;s,y)=(t-s)^{\beta-\gamma}\varrho^\beta_{\gamma-\beta}(t,x;s,y)
\preceq|x-x'|^{\beta-\gamma}\varrho^\beta_{\gamma-\beta}(t,x;s,y),
$$
and also
$$
|q_0(t,x';s,y)|\preceq|x-x'|^{\beta-\gamma}\varrho^\beta_{\gamma-\beta}(t,x';s,y).
$$
Suppose now that
\begin{align}
|x-x'|\leq|t-s|.\label{EE1}
\end{align}
We can write
\begin{align*}
|q_0(t,x;s,y)-q_0(t,x';s,y)|&\leq|a(t,x)-a(t,y)|\cdot|\Delta^{\frac{1}{2}}_xp_0(t,x;s,y)
-\Delta^{\frac{1}{2}}_{x'}p_0(t,x';s,y)|\\
&+|a(t,x)-a(t,x')|\cdot|\Delta^{\frac{1}{2}}_{x'}p_0(t,x';s,y)|\\
&+|b(t,x)-b(t,y)|\cdot|\nabla_xp_0(t,x;s,y)-\nabla_{x'}p_0(t,x';s,y)|\\
&+|b(t,x)-b(t,x')|\cdot|\nabla_{x'}p_0(t,x';s,y)|\\
&=:I_1+I_2+I_3+I_4.
\end{align*}
For $I_1$, by (\ref{EU6}) and the mean value theorem, we have for some $\theta\in[0,1]$,
$$
I_1\preceq\big\{|x-y|^\beta\wedge1\}|x-x'|(|x+\theta(x'-x)-y|+|t-s|)^{-d-2}.
$$
By (\ref{EE1}), we have
$$
|x-y|+|t-s|\leq|x+\theta(x'-x)-y|+2|t-s|.
$$
Hence,
\begin{align*}
I_1&\preceq\big\{|x-y|^\beta\wedge1\}|x-x'|(|x-y|+|t-s|)^{-d-2}\\
&\preceq|x-x'|^{\beta-\gamma}\frac{|t-s|^{1+\gamma-\beta}\big\{|x-y|^\beta\wedge1\}}{|x-y|+|t-s|}(|x-y|+|t-s|)^{-d-1}\\
&\preceq|x-x'|^{\beta-\gamma}|t-s|^\gamma(|x-y|+|t-s|)^{-d-1}
=|x-x'|^{\beta-\gamma}\varrho^0_{\gamma}(t,x;s,y).
\end{align*}
By (\ref{EU4}), we have
$$
I_2\preceq |x-x'|^\beta(|x'-y|+|t-s|)^{-d-1}\preceq
|x-x'|^{\beta-\gamma}\varrho^0_{\gamma}(t,x';s,y).
$$
Similarly, we have
$$
I_3\preceq|x-x'|^{\beta-\gamma}\varrho^0_{\gamma}(t,x;s,y).
$$
$$
I_4\preceq|x-x'|^{\beta-\gamma}\varrho^0_{\gamma}(t,x';s,y).
$$
Combining the above calculations, we obtain (\ref{ES5}).

Now, by definition (\ref{ES6}), (\ref{ES5}) and Lemma \ref{Le3}, we have for $n\in\mN$,
\begin{align*}
&|q_n(t,x;s,y)-q_n(t,x';s,y)|\preceq\int_s^t\!\!\!\int_{\mR^d}|q_0(t,x;r,z)-q_0(t,x';r,z)|q_{n-1}(r,z;s,y)\dif z\dif r\\
&\qquad\preceq\frac{(C_d\Gamma(\beta))^{n}}{\Gamma(n\beta)}\Big(|x-x'|^{\beta-\gamma}\wedge 1\Big)\int_s^t\!\!\!\int_{\mR^d}
\Big\{(\varrho_{\gamma}^0+\varrho_{\gamma-\beta}^{\beta})(t,x;r,z)+(\varrho_{\gamma}^0+\varrho_{\gamma-\beta}^{\beta})(t,x';r,z)\Big\}\\
&\qquad\qquad\qquad\qquad\qquad\qquad\qquad\times\Big\{\varrho_{n\beta}^{0}(r,z;s,y)+\varrho_{(n-1)\beta}^{\beta}(r,z;s,y)\Big\}\dif z\dif r,
\end{align*}
which yields the result by Lemma \ref{Le2}.
\end{proof}

Basing on the above two lemmas, we have
\bt\label{T3.4}
The function $q(t,x;s,y):=\sum_{n=0}^{\infty}q_n(t,x;s,y)$ solves
the integral equation (\ref{Eq}). Moreover, $q(t,x;s,y)$ has the following estimates:
\begin{align}
|q(t,x;s,y)|\preceq\varrho^\beta_0(t,x;s,y)+\varrho^0_\beta(t,x;s,y), \label{eq3}
\end{align}
and for any $\gamma\in(0,\beta)$,
\begin{align}
&|q(t,x;s,y)-q(t,x';s,y)|\no\\
&\qquad\preceq\Big(|x-x'|^{\beta-\gamma}\wedge 1\Big)\Big\{(\varrho^0_\gamma+\varrho^\beta_{\gamma-\beta})(t,x;s,y)
+(\varrho^0_\gamma+\varrho^\beta_{\gamma-\beta})(t,x';s,y)\Big\}.
\label{eq4}
\end{align}
In the case of $a(t,x)=a(t)$ independent of $x$, we have
\begin{align}
|q(t,x;s,y)|\preceq \varrho^0_\beta(t,x;s,y). \label{eq5}
\end{align}
\et
\begin{proof}
By Lemma \ref{Le3}, one sees that
\begin{align*}
\sum_{n=0}^{\infty}|q_n(t,x;s,y)|&\leq\sum_{n=0}^\infty
\frac{(C_d\Gamma(\beta))^{n+1}}{\Gamma((n+1)\beta)}\left(\varrho^0_{(n+1)\beta}(t,x;s,y)
+\varrho^\beta_{n\beta}(t,x;s,y)\right)\\
&\leq\left\{\sum_{n=0}^\infty\frac{(C_d\Gamma(\beta))^{n+1}}{\Gamma((n+1)\beta)}\right\}
\left(\varrho^0_{\beta}(t,x;s,y)+\varrho^\beta_0(t,x;s,y)\right).
\end{align*}
Since the series is convergent, we obtain (\ref{eq3}). Similarly, estimate (\ref{eq4}) follows by Lemma \ref{Le4}.
Moreover, by (\ref{ES6}) we have
$$
\sum_{n=0}^{m+1}q_n(t,x;s,y)=q_0(t,x;s,y)+\int^t_s\!\!\!\int_{\mR^d}q_0(t,x;r,z)\sum_{n=0}^mq_n(r,z;s,y)\dif z\dif r,
$$
which yields (\ref{Eq}) by taking limits $m\to\infty$ for both sides.

In the case of $a(t,x)=a(t)$, we use (\ref{EW3}) to repeat the above proof, and obtain (\ref{eq5}).
\end{proof}

For $r\in(s,t)$, let us set
$$
\phi(t,x,r):=\phi_{s,y}(t,x,r):=\int_{\mR^d}p_0(t,x;r,z)q(r,z;s,y)\dif z,
$$
and
$$
\varphi(t,x):=\int^t_s\phi(t,x,r)\dif r=\int^t_s\!\!\!\int_{\mR^d}p_0(t,x;r,z)q(r,z;s,y)\dif z\dif r,
$$
where the integral is taken in the generalized sense, i.e.,
$$
\phi(t,x,r)=\lim_{\eps\downarrow 0}\int_{|x-z|\geq\eps}p_0(t,x;r,z)q(r,z;s,y)\dif z.
$$
Notice that by (\ref{E1}), (\ref{eq3})  and (\ref{EU7}),
\begin{align}
|\phi_{s,y}(t,x,r)|&\leq\int_{\mR^d}p_0(t,x;r,z)|q(r,z;s,y)|\dif z
\preceq\int_{\mR^d}\varrho^0_1(t,x;r,z)(\varrho^0_\beta+\varrho^\beta_0)(r,z;s,y)\dif z\no\\
&\preceq \Big((t-r)^{\beta}+(r-s)^{\beta}+(t-r)(r-s)^{\beta -1}\Big)\varrho^0_0(t,x;s,y)+\varrho^\beta_0(t,x;s,y).\label{eq7}
\end{align}

Below we study the smoothness of $(t,x)\mapsto\varphi(t,x)$.
\bl
For all $x\not=y\in\mR^d$ and  almost all $t>s$, we have
\begin{align}
\p_t \varphi(t,x)=q(t,x;s,y)+\int^t_s\!\!\!\int_{\mR^d}\sL_{a,b}^x(t,z)p_0(t,\cdot;r,z)(x)q(r,z;s,y)\dif z\dif r.\label{ES8}
\end{align}
\el
\begin{proof}
(Claim 1): For $r\in(s,t)$, we have
\begin{align}
\p_t \phi_{s,y}(t,x,r)=\int_{\mR^d}\p_tp_0(t,x;r,z)q(r,z;s,y)\dif z.\label{EW1}
\end{align}
{\it Proof of Claim 1}: Write
\begin{align*}
\frac{\phi_{s,y}(t+\eps,x,r)-\phi_{s,y}(t,x,r)}{\eps}&=\frac{1}{\eps}\int_{\mR^d}\big(p_0(t+\eps,x;r,z)-p_0(t,x;r,z)\big)q(r,z;s,y)\dif z\no\\
&=\int_{\mR^d}\!\Bigg(\!\int_0^1 \p_t p_0(t+\theta \eps,x;r,z)\dif \theta \Bigg)q(r,z;s,y)\dif z.
\end{align*}
By (\ref{EU5}) and (\ref{EU4}), we have  for $|\eps|<\frac{t-r}{2}$,
$$
|\p_t p_0(t+\theta \eps,x;r,z)|\preceq (|x-z|+t+\theta\eps-r)^{-d-1}\preceq (|x-z|+(t-r))^{-d-1},
$$
which together with (\ref{eq3}) yields
\begin{align*}
|\p_t p_0(t+\theta \eps,x;r,z)q(r,z;s,y)|
\preceq \varrho_0^0(t,x;r,z)(\varrho^0_\beta+\varrho^\beta_0)(r,z;s,y)=:g(z).
\end{align*}
By (\ref{EU7}), one sees that
$$
\int_{\mR^d}g(z)\dif z<+\infty.
$$
Hence, by the dominated convergence theorem, we have
\begin{align*}
\lim_{\eps\to 0}\frac{\phi_{s,y}(t+\eps,x,r)-\phi_{s,y}(t,x,r)}{\eps}
=\int_{\mR^d}\p_tp_0(t,x;r,z)q(r,z;s,y)\dif z,
\end{align*}
and (\ref{EW1}) is proven.

\vspace{5mm}

(Claim 2): For $x\neq y$, we have
\begin{align}
\int^t_s\!\!\!\int^{t'}_s|\p_{t'} \phi_{s,y}(t',x,r)|\dif r\dif t'<+\infty.\label{EU9}
\end{align}
{\it Proof of Claim 2}:
By (\ref{EW1}), we have
\begin{align}
|\p_{t'} \phi_{s,y}(t',x,r)|&\leq\int_{\mR^d}|\p_{t'}p_0(t',x;r,z)|\cdot|q(r,z;s,y)-q(r,x;s,y)|\dif z\no\\
&\quad+|q(r,x;s,y)|\left|\int_{\mR^d}\p_{t'}p_0(t',x;r,z)\dif z\right|\no\\
&=:Q^{(1)}_{s,y}(t',x,r)+Q^{(2)}_{s,y}(t',x,r).\label{EP4}
\end{align}
For $Q^{(1)}_{s,y}(t',x,r)$, by (\ref{EU404}) and (\ref{eq4}), we have
\begin{align}
\int^t_s\!\!\!\int^{t'}_sQ^{(1)}_{s,y}(t',x,r)\dif r\dif t'&\preceq
\int^t_s\!\!\!\int^{t'}_s\!\!\!\int_{\mR^d}\varrho^{\beta-\gamma}_0(t',x;r,z)
(\varrho^0_\gamma+\varrho^\beta_{\gamma-\beta})(r,x;s,y)\dif z\dif r\dif t'\no\\
&\quad+\int^t_s\!\!\!\int^{t'}_s\!\!\!\int_{\mR^d}\varrho^{\beta-\gamma}_0(t',x;r,z)
(\varrho^0_\gamma+\varrho^\beta_{\gamma-\beta})(r,z;s,y)\dif z\dif r\dif t'\no\\
&\preceq \int^t_s\!\!\!\int^{t'}_s(t'-r)^{\beta-\gamma-1}(\varrho^0_\gamma+\varrho^\beta_{\gamma-\beta})(r,x;s,y)\dif r\dif t'\no\\
&\quad+\int^t_s(\varrho^0_\beta+\varrho^\beta_0+\varrho^{\beta-\gamma}_\gamma)(t',x;s,y)\dif t'\no\\
&\preceq \frac{1}{|x-y|^{d+1}}\int^t_s\!\!\!\int^{t'}_s(t'-r)^{\beta-\gamma-1}((r-s)^\gamma+(r-s)^{\gamma-\beta})\dif r\dif t'\no\\
&\quad+\frac{1}{|x-y|^{d+1}}\int^t_s((t'-s)^\gamma+1+(t'-s)^\beta)\dif t'<+\infty.\label{EP5}
\end{align}
For $Q^{(2)}_{s,y}(t',x,r)$, by (\ref{EW6}) and (\ref{eq3}) we have
\begin{align}
\int^t_s\!\!\!\int^{t'}_sQ^{(2)}_{s,y}(t',x,r)\dif r\dif t'&\preceq\int^t_s\!\!\!\int^{t'}_s
(\varrho^0_\beta+\varrho^\beta_0)(r,x;s,y)(t'-r)^{\beta-1}\dif r\dif t'<+\infty.\label{EP6}
\end{align}
Combining (\ref{EP4})-(\ref{EP6}), we obtain (\ref{EU9}).

\vspace{5mm}

(Claim 3): For fixed $r,x,s,y$, we have
\begin{align}
\lim_{t\downarrow r}\phi_{s,y}(t,x,r)=q(r,x;s,y).\label{EW11}
\end{align}
{\it Proof of Claim 3}:
By (\ref{EW66}), it suffices to prove that
$$
\lim_{t\downarrow r}\Bigg|\int_{\mR^d}p_0(t,x;r,z)(q(r,z;s,y)- q(r,x;s,y))\dif z \Bigg|=0.
$$
Notice that for any $\delta>0$,
\begin{align*}
&\Bigg|\int_{\mR^d}p_0(t,x;r,z)(q(r,z;s,y)- q(r,x;s,y))\dif z \Bigg|\\
&\quad \leq\int_{|x-z|\leq\delta}p_0(t,x;r,z)|q(r,z;s,y)-q(r,x;s,y)|\dif z\\
&\quad +\int_{|x-z|>\delta}p_0(t,x;r,z)|q(r,z;s,y)-q(r,x;s,y)|\dif z\\
&\quad =: J_1(\delta,t,r)+J_2(\delta,t,r).
\end{align*}
For any $\eps>0$, by (\ref{eq4}), there exists a $\delta=\delta(r,x,s,y)>0$ such that
for all $|x-z|\leq\delta$,
$$
|q(r,z;s,y)-q(r,x;s,y)|\leq \eps.
$$
Thus,
$$
J_1(\delta,t,r)\leq \varepsilon\int_{|x-z|\leq \delta}p_0(t,x;r,z)\dif z
\leq\varepsilon\int_{\mR^d}p_0(t,x;r,z)\dif z\preceq\eps\int_{\mR^d}\varrho^0_1(t,x;r,z)\dif z\leq\eps.
$$
On the other hand, we have
\begin{align*}
J_2(\delta,t,r)&\stackrel{(\ref{E1})}{\preceq}(t-r)\int_{|x-z|>\delta}\frac{|q(r,z;s,y)|+|q(r,x;s,y)|}{|x-z|^{d+1}}\dif z\\
&\leq(t-r)\left(\delta^{-d-1}\int_{\mR^d}|q(r,z;s,y)|\dif z
+|q(r,x;s,y)|\int_{|z|>\delta}|z|^{-d-1}\dif z\right),
\end{align*}
which, by (\ref{eq3}) and (\ref{ES4}), converges to zero as $t\downarrow r$. The claim (\ref{EW11}) is thus proved.

\vspace{5mm}

Now, by integration by parts formula and (\ref{EW11}), we have
$$
\int^t_r\p_{r'} \phi_{s,y}(r',x,r)\dif r'=\phi_{s,y}(t,x,r)-q(r,x;s,y).
$$
Integrating both sides with respect to $r$ from $s$ to $t$, and then by (\ref{EU9}) and Fubini's theorem, we obtain
\begin{align*}
\varphi(t,x)-\int^t_sq(r,x;s,y)\dif r&=\int^t_s\!\!\!\int^t_r\p_{r'} \phi_{s,y}(r',x,r)\dif r'\dif r
\stackrel{(\ref{EU9})}{=}\int^t_s\!\!\!\int^{r'}_s\p_{r'} \phi_{s,y}(r',x,r)\dif r\dif r'\\
&\stackrel{(\ref{EW1}), (\ref{ES2})}{=}\int^t_s\!\!\!\int^{r'}_s\!\!\!\int_{\mR^d}\sL_{a,b}^x(r',z)
p_0(r',\cdot;r,z)(x)q(r,z;s,y)\dif z\dif r\dif r',
\end{align*}
which in turn implies (\ref{ES8}) by the Lebesgue differential theorem.
\end{proof}
\bl\label{Le6}
For all $t>s$ and $x\not=y$, we have
\begin{align}
\nabla_x \varphi(t,x)&=\int^t_s\!\!\!\int_{\mR^d}\nabla_xp_0(t,x;r,z)q(r,z;s,y)\dif z\dif r,\label{EP2}\\
\Delta^{\frac{1}{2}}_x \varphi(t,x)&=\int^t_s\!\!\!\int_{\mR^d}\Delta^{\frac{1}{2}}_x
p_0(t,x;r,z)q(r,z;s,y)\dif z\dif r,\label{EP1}
\end{align}
where the integrals are understood in the sense of iterated integrals.
\el
\begin{proof}
First of all, for fixed $s<r<t$, since
$$
(x,z)\mapsto p_0(t,x;r,z)\in C^\infty_b(\mR^d\times\mR^d)
$$
and
$$
z\mapsto q(r,z;s,y)\in C_b(\mR^d),
$$
by Lemma \ref{Le1}, it is easy to see that
\begin{align}
\nabla_x\phi_{s,y}(t,x,r)=\int_{\mR^d}\nabla_xp_0(t,x;r,z)q(r,z;s,y)\dif z,\label{EY2}
\end{align}
and
\begin{align}
\Delta^{\frac{1}{2}}_x\phi_{s,y}(t,x,r)=\int_{\mR^d}\Delta^{\frac{1}{2}}_xp_0(t,x;r,z)q(r,z;s,y)\dif z.\label{EY3}
\end{align}

(1) We prove the following claim: For any $t>s$ and $x\not=y$, there exists a $p>1$ such that
\begin{align}
\sup_{|w|\leq|x-y|/2}I(p,w)<\infty,\mbox{ where } I(p,w):=\int^t_s|\nabla_x\phi_{s,y}(t,x+w;r)|^p\dif r.\label{EY1}
\end{align}
{\it Proof of Claim}: By (\ref{EY2}) and (\ref{EW7}) we have
\begin{align*}
I(p,w)&\preceq\int^t_s\left|\int_{\mR^d}\nabla_x p_0(t,x+w;r,z)(q(r,z;s,y)-q(r,x+w;s,y))\dif z\right|^p\dif r\\
&\quad+\int^t_s\left|\int_{\mR^d}\nabla_x p_0(t,x+w;r,z)\dif z\right|^p|q(r,x+w;s,y)|^p\dif r\\
&\stackrel{(\ref{EU4}) (\ref{eq4})}{\preceq}\int^t_s\left(\int_{\mR^d}\varrho^{\beta-\gamma}_0(t,x+w;r,z)
(\varrho^0_\gamma+\varrho^\beta_{\gamma-\beta})(r,z;s,y)\dif z\right)^p\dif r\\
&\qquad+\int^t_s\left(\int_{\mR^d}\varrho^{\beta-\gamma}_0(t,x+w;r,z)
(\varrho^0_\gamma+\varrho^\beta_{\gamma-\beta})(r,x+w;s,y)\dif z\right)^p\dif r\\
&\stackrel{(\ref{EW7})}{}+\int^t_s(t-r)^{p(\beta-1)}(\varrho^0_\beta+\varrho^\beta_0)^p(r,x+w;s,y)\dif r\\
&=:I_1(p,w)+I_2(p,w)+I_3(p,w).
\end{align*}
For $I_1(p,w)$, it follows by (\ref{eq31}) that for some $p>1$,
$$
\sup_{|w|\leq|x-y|/2}I_1(p,w)<+\infty.
$$
For $I_2(p,w)$, by definition (\ref{Def2}) and (\ref{ES4}), we have for all $|w|\leq|x-y|/2$,
\begin{align*}
I_2(p,w)&\preceq\int^t_s\left(\int_{\mR^d}\varrho^{\beta-\gamma}_0(t,x+w;r,z)\dif z\right)^p
\left(\frac{(r-s)^\gamma}{|x+w-y|^{d+1}}+\frac{(r-s)^{\gamma-\beta}}{|x+w-y|^{d+1}}\right)^p\dif r\\
&\preceq\int^t_s(t-r)^{p(\beta-\gamma-1)}
\left(\frac{1}{|x-y|^{d+1}}+\frac{(r-s)^{\gamma-\beta}}{|x-y|^{d+1}}\right)^p\dif r<+\infty,
\end{align*}
provided $p<\frac{1}{1+\gamma-\beta}\wedge\frac{1}{\beta-\gamma}$.
Similarly, for $p<\frac{1}{1-\beta}$, we have
\begin{align}
\sup_{|w|\leq|x-y|/2}I_3(p,w)<+\infty.
\end{align}
Thus, we obtain (\ref{EY1}).

\vspace{5mm}

Now, for any $e_i=(0,\cdots,1,\cdots,0)\in\mR^d$, we can write
\begin{align*}
\frac{\varphi(t,x+\eps e_i)-\varphi(t,x)}{\eps}
=\int^t_s\!\!\!\int^1_0\p_{x_i}\phi(t,x+\theta\eps e_i,r)\dif\theta\dif r.
\end{align*}
By (\ref{EY1}) one can take limits to get
\begin{align*}
\p_{x_i}\varphi(t,x)=
\lim_{\eps\to 0}\frac{\varphi(t,x+\eps e_i)-\varphi(t,x)}{\eps}
&=\int^t_s\!\!\!\int^1_0\lim_{\eps\to 0}\p_{x_i}\phi(t,x+\theta\eps e_i,r)\dif\theta\dif r
=\int^t_s\p_{x_i}\phi(t,x,r)\dif r,
\end{align*}
and (\ref{EP2}) is proven.

(2) Next, we prove (\ref{EP1}). Recalling the definition of $\phi_{s,y}$, we have
\begin{align*}
&\nabla_x\phi_{s,y}(t,x+w,r)-\nabla_x\phi_{s,y}(t,x,r)\\
&\qquad=\int_{\mR^d}(\nabla_xp_0(t,x+w;r,z)-\nabla_xp_0(t,x;r,z))q(r,z;s,y)\dif z\\
&\qquad=\int_{\mR^d}\Big(\nabla_xp_0(t,x+w;r,z)(q(r,z;s,y)-q(r,x+w;s,y))\\
&\qquad\qquad-\nabla_xp_0(t,x;r,z)(q(r,z;s,y)-q(r,x;s,y))\Big)\dif z\\
&\qquad\quad+\Bigg\{q(r,x+w;s,y)\int_{\mR^d}\nabla_xp_0(t,x+w;r,z)\dif z\\
&\qquad\qquad-q(r,x;s,y)\int_{\mR^d}\nabla_xp_0(t,x;r,z)\dif z\Bigg\}\\
&\qquad=:\int_{\mR^d}Q(t,x;r,z;s,y;w)\dif z+R(r,t,x,w;s,y).
\end{align*}
We now prove the following claim: For any $\gamma\in(0,\beta)$ and $\sigma\in(0,\beta-\gamma)$,
\begin{align}
&|Q(t,x;r,z;s,y;w)|\no\\
&\quad\preceq|w|^\sigma\varrho^{\beta-\gamma}_{-\sigma}(t,x+w;r,z)
\Big((\varrho^\beta_{\gamma-\beta}+\varrho^0_\gamma)(r,x+w;s,y)+
(\varrho^\beta_{\gamma-\beta}+\varrho^0_\gamma)(r,z;s,y)\Big)\no\\
&\quad+|w|^\sigma\varrho^{\beta-\gamma}_{-\sigma}(t,x;r,z)
\Big((\varrho^\beta_{\gamma-\beta}+\varrho^0_\gamma)(r,x;s,y)+
(\varrho^\beta_{\gamma-\beta}+\varrho^0_\gamma)(r,z;s,y)\Big)\no\\
&\quad+|w|^\sigma\varrho^0_{\beta-\gamma-\sigma}(t,x;r,z)
\Big((\varrho^\beta_{\gamma-\beta}+\varrho^0_\gamma)(r,x+w;s,y)+
(\varrho^\beta_{\gamma-\beta}+\varrho^0_\gamma)(r,z;s,y)\Big)\no\\
&\quad+|w|^\sigma\varrho^0_{\beta-\gamma-\sigma}(t,x;r,z)
\Big((\varrho^\beta_{\gamma-\beta}+\varrho^0_\gamma)(r,x+w;s,y)+
(\varrho^\beta_{\gamma-\beta}+\varrho^0_\gamma)(r,x;s,y)\Big),\label{EY4}
\end{align}
and for $w\in\mR^d$,
\begin{align}
R(r,t,x;s,y;w)&\preceq |w|^{\beta-\gamma}(t-r)^{\beta-1}\Big\{(\varrho^0_\gamma+\varrho^\beta_{\gamma-\beta})(r,x+w;s,y)
+(\varrho^0_\gamma+\varrho^\beta_{\gamma-\beta})(r,x;s,y)\Big\}\no\\
&\quad+|w|^\gamma(t-r)^{\beta-\gamma-1}(\varrho^0_\beta+\varrho^{\beta}_0)(r,x;s,y).\label{EQ4}
\end{align}
{\it Proof of Claim}:
Case 1: We assume
$$
|w|>|t-r|.
$$
By (\ref{eq4}) we have
\begin{align*}
&|\nabla_xp_0(t,x;r,z)(q(r,z;s,y)-q(r,x;s,y))|\\
&\quad\preceq\varrho^0_0(t,x;r,z)
(|x-z|^{\beta-\gamma}\wedge 1)\Big((\varrho^\beta_{\gamma-\beta}+\varrho^0_\gamma)(r,x;s,y)+
(\varrho^\beta_{\gamma-\beta}+\varrho^0_\gamma)(r,z;s,y)\Big)\\
&\quad\preceq |w|^\sigma\varrho^{\beta-\gamma}_{-\sigma}(t,x;r,z)
\Big((\varrho^\beta_{\gamma-\beta}+\varrho^0_\gamma)(r,x;s,y)+
(\varrho^\beta_{\gamma-\beta}+\varrho^0_\gamma)(r,z;s,y)\Big),
\end{align*}
and also
\begin{align*}
&|\nabla_xp_0(t,x+w;r,z)(q(r,z;s,y)-q(r,x+w;s,y))|\\
&\quad\preceq |w|^\sigma\varrho^{\beta-\gamma}_{-\sigma}(t,x+w;r,z)
\Big((\varrho^\beta_{\gamma-\beta}+\varrho^0_\gamma)(r,x+w;s,y)+
(\varrho^\beta_{\gamma-\beta}+\varrho^0_\gamma)(r,z;s,y)\Big).
\end{align*}
Case 2: We assume
$$
|w|\leq|t-r|.
$$
Noticing that
$$
|x+w-z|\leq|x-z|+|w|\leq|x-z|+|t-r|
$$
and
$$
|x-z|\leq|x+w-z|+|w|\leq|x+w-z|+|t-r|,
$$
we have for any $\theta_0\in(0,1)$,
\begin{align*}
|w|\cdot|\nabla^2_xp_0(t,x+\theta_0 w;r,z)|\cdot|x+w-z|^{\beta-\gamma}\preceq
|w|^\sigma\varrho^0_{\beta-\gamma-\sigma}(t,x;r,z).
\end{align*}
Hence, for some $\theta_0\in(0,1)$,
\begin{align*}
&|(\nabla_xp_0(t,x+w;r,z)-\nabla_xp_0(t,x;r,z))(q(r,z;s,y)-q(r,x+w;s,y))|\\
&\qquad\preceq |w|\cdot|\nabla^2_xp_0(t,x+\theta_0 w;r,z)|\cdot|x+w-z|^{\beta-\gamma}\\
&\qquad\times\Big((\varrho^\beta_{\gamma-\beta}+\varrho^0_\gamma)(r,x+w;s,y)+
(\varrho^\beta_{\gamma-\beta}+\varrho^0_\gamma)(r,z;s,y)\Big)\\
&\qquad\preceq|w|^\sigma\varrho^0_{\beta-\gamma-\sigma}(t,x;r,z)\Big((\varrho^\beta_{\gamma-\beta}+\varrho^0_\gamma)(r,x+w;s,y)+
(\varrho^\beta_{\gamma-\beta}+\varrho^0_\gamma)(r,z;s,y)\Big).
\end{align*}
Similarly, we have
\begin{align*}
&|\nabla_xp_0(t,x;r,z)(q(r,x;s,y)-q(r,x+w;s,y))|\\
&\quad\preceq|w|^\sigma\varrho^0_{\beta-\gamma-\sigma}(t,x;r,z)
\Big((\varrho^\beta_{\gamma-\beta}+\varrho^0_\gamma)(r,x+w;s,y)+
(\varrho^\beta_{\gamma-\beta}+\varrho^0_\gamma)(r,x;s,y)\Big).
\end{align*}
Combining the above calculations, we obtain (\ref{EY4}).
As for (\ref{EQ4}), it follows by Lemma \ref{Le1} and Theorem \ref{T3.4}.

\vspace{5mm}

(3) Now, we are ready to prove the following claim: For any $t>s$ and $x\not=y$, there exists a $p>1$ such that
\begin{align}
\sup_{\eps\leq|x-y|/2}J(p,\eps)<+\infty,\mbox{ where }J(p,\eps):=\int^t_s\left|\int_{|w|\geq\eps}
\frac{\phi_{s,y}(t,x+w,r)-\phi_{s,y}(t,x,r)}{|w|^{d+1}}\dif w\right|^p\dif r.\label{EY0}
\end{align}
{\it Proof of Claim}: Notice that
\begin{align*}
J(p,\eps)&\preceq\int^t_s\left|\int_{\eps\leq|w|\leq|x-y|/2}\frac{\phi_{s,y}(t,x+w,r)-\phi_{s,y}(t,x,r)
-w\cdot\nabla_x\phi_{s,y}(t,x,r)}{|w|^{d+1}}\dif w\right|^p\dif r\\
&\quad+\int^t_s\left|\int_{|w|>|x-y|/2}\frac{\phi_{s,y}(t,x+w,r)-\phi_{s,y}(t,x,r)}{|w|^{d+1}}\dif w\right|^p\dif r
=:J_1(p,\eps)+J_2(p).
\end{align*}
For $J_1(p,\eps)$, observe that
\begin{align*}
J_1(p,\eps)&=\int^t_s\left|\int_{\eps\leq|w|\leq|x-y|/2}\frac{w}{|w|^{d+1}}\cdot
\left(\int^1_0\Big(\nabla_x\phi_{s,y}(t,x+\theta w,r)-\nabla_x\phi_{s,y}(t,x,r)\Big)\dif\theta\right)\dif w\right|^p\dif r\\
&\preceq\int^t_s\left(\int_{\eps\leq|w|\leq|x-y|/2}
\!\int^1_0\frac{\int_{\mR^d}|Q(t,x;r,z;s,y;\theta w)|\dif z}{|w|^d}\dif\theta\dif w\right)^p\dif r\\
&\qquad+\int^t_s\left(\int_{\eps\leq|w|\leq|x-y|/2}
\!\int^1_0\frac{R(r,t,x;s,y;\theta w)}{|w|^d}\dif\theta\dif w\right)^p\dif r.
\end{align*}
Using (\ref{EY4}), (\ref{EQ4}) and by (\ref{eq31}), as in proving (\ref{EY1}), one has that for some $p>1$,
$$
\sup_{\eps\leq|x-y|/2}J_1(p,\eps)<+\infty.
$$
For $J_2(p)$, we have by (\ref{eq7}) that for some $p>1$,
$$
J_2(p)\leq\int^t_s\left|\int_{|w|>|x-y|/2}\frac{|\phi_{s,y}(t,x+w,r)|+|\phi_{s,y}(t,x,r)|}{|w|^{d+1}}\dif w\right|^p\dif r<+\infty.
$$
Thus, (\ref{EY0}) is proven.

\vspace{5mm}

(4) Now by (\ref{EY0}), one has
\begin{align*}
\Delta^{\frac{1}{2}}_x\varphi(t,x)&=\lim_{\eps\downarrow 0}\int_{|w|\geq\eps}\!\int^t_s
\frac{\phi_{s,y}(t,x+w,r)-\phi_{s,y}(t,x,r)}{|w|^{d+1}}\dif r\dif w\\
&=\lim_{\eps\downarrow 0}\int^t_s\!\!\!\int_{|w|\geq\eps}
\frac{\phi_{s,y}(t,x+w,r)-\phi_{s,y}(t,x,r)}{|w|^{d+1}}\dif w\dif r\\
&=\int^t_s\lim_{\eps\downarrow 0}\int_{|w|\geq\eps}
\frac{\phi_{s,y}(t,x+w,r)-\phi_{s,y}(t,x,r)}{|w|^{d+1}}\dif w\dif r
=\int^t_s\Delta^{\frac{1}{2}}_x\phi_{s,y}(t,x,r)\dif r,
\end{align*}
which together with (\ref{EY3}) yields (\ref{EP1}).
\end{proof}

Now we prove the following main result of this section.
\bt\label{T1}
Assume that $a,b\in\mH^\beta$ for some $\beta \in (0,\frac{1}{4}]$ and satisfy (\ref{Con}). Then there exists a unique nonnegative continuous
function $p_{a,b}(t,x;s,y)$ with
\begin{align}
\int_{\mR^d}p_{a,b}(t,x;s,y)\dif y=1,\label{EK1}
\end{align}
and satisfying that
\begin{enumerate}[(i)]
\item For all $x\not=y\in\mR^d$ and almost all $t>s$,
\begin{align}
\p_tp_{a,b}(t,x;s,y)=\sL^x_{a,b}(t,x)p_{a,b}(t,\cdot;s,y)(x). \label{eq15}
\end{align}
\item
For any bounded continuous function $f$,
\begin{align}\label{eq14}
\lim_{t\downarrow s}\int_{\mR^d}p_{a,b}(t,x;s,y)f(y)\dif y=f(x),
\end{align}
and
\begin{align}\label{eq141}
\lim_{t\downarrow s}\int_{\mR^d}p_{a,b}(t,x;s,y)f(x)\dif x=f(y).
\end{align}
\item  For all $0\leq s<t\leq 1$ and $x,y\in\mR^d$,
\begin{align}
p_{a,b}(t,x;s,y)\preceq \varrho^0_1(t,x;s,y); \label{eq16}
\end{align}
and in the case of $a(t,x)=a(t)$ independent of $x$, we have
\begin{align}
p_{a,b}(t,x;s,y)\succeq\varrho^0_1(t,x;s,y). \label{eq18}
\end{align}
\item For any $\gamma\in(0,1)$,
\begin{align}\label{eq19}
|p_{a,b}(t,x;s,y)-p_{a,b}(t,x';s,y)|\preceq (|x-x'|^{\gamma}\wedge 1)\Big\{\varrho_{1-\gamma}^0(t,x;s,y)+\varrho_{1-\gamma}^0(t,x';s,y)\Big\},
\end{align}
and
\begin{align}
|\nabla_xp_{a,b}(t,x;s,y)|+|\Delta^{\frac{1}{2}}_xp_{a,b}(t,x;s,y)|\preceq \varrho^0_0(t,x;s,y).\label{eq17}
\end{align}
\item If $a,b\in C\big([0,\infty);L_{\mathrm{loc}}^{1}(\mR^d)\big)$, then for all $f,g\in C^\infty_0(\mR^d)$,
\begin{align}\label{eq233}
\lim_{t\downarrow s}\frac{1}{t-s}\int_{\mR^d}g(x)\Big(P^{a,b}_{t,s}f(x)-f(x)\Big)\dif x
=\int_{\mR^d}g(x)\sL^x_{a,b}(s,x)f(x)\dif x,\ \ s\geq 0,
\end{align}
where
$$
P^{a,b}_{t,s}f(x):=\int_{\mR^d}p_{a,b}(t,x;s,y)f(y)\dif y.
$$
\end{enumerate}
\et
\begin{proof}
(i) First, we prove (\ref{eq15}). By equation (\ref{eq2}), we have that for all $x,y\in\mR^d$ and almost all $t>s$,
\begin{align*}
\p_tp_{a,b}(t,x;s,y)&\stackrel{(\ref{ES8})}{=}\p_tp_0(t,x;s,y)
+q(t,x;s,y)+\int^t_s\!\!\!\int_{\mR^d}\sL_{a,b}^x(t,y)
p_0(t,\cdot;r,z)(x)q(r,z;s,y)\dif z\dif r\\
&\stackrel{(\ref{ES2})}{=}\sL^x_{a,b}(t,y)p_0(t,\cdot;s,y)(x)+q_0(t,x;s,y)
+\int^t_s\!\!\!\int_{\mR^d}q_0(t,x;r,z)q(r,z;s,y)\dif z\dif r\\
&\quad+\int^t_s\!\!\!\int_{\mR^d}\sL_{a,b}^x(t,y)
p_0(t,\cdot;r,z)(x)q(r,z;s,y)\dif z\dif r.
\end{align*}
Recalling that
\begin{align}
q_0(t,x;s,y)=\left(\sL^x_{a,b}(t,x)-\sL^x_{a,b}(t,y)\right)p_0(t,\cdot;s,y)(x),\label{EW4}
\end{align}
we further have
\begin{align*}
\p_tp_{a,b}(t,x;s,y)=\sL^x_{a,b}(t,x)p_0(t,\cdot;s,y)(x)
+\int^t_s\!\!\!\int_{\mR^d}\sL_{a,b}^x(t,x)p_0(t,\cdot;r,z)(x)q(r,z;s,y)\dif z\dif r,
\end{align*}
which together with (\ref{EP2}) and (\ref{EP1}) yields (\ref{eq15}).

(ii) We prove (\ref{eq14}) and (\ref{eq141}). As in proving (\ref{EW11}), we can prove that for any bounded continuous function $f$,
$$
\lim_{t\downarrow s}\int_{\mR^d}p_0(t,x;s,y)f(y)\dif y=f(x),
$$
and
$$
\lim_{t\downarrow s}\int_{\mR^d}p_0(t,x;s,y)f(x)\dif x=f(y).
$$
Moreover, by (\ref{eq7}) we also have
\begin{align*}
&\left|\int_{\mR^d}\!\int_s^t\!\!\!\int_{\mR^d}p_0(t,x;r,z)q(r,z;s,y)f(y)\dif z\dif r\dif y\right|
+\left|\int_{\mR^d}\!\int_s^t\!\!\!\int_{\mR^d}p_0(t,x;r,z)q(r,z;s,y)f(x)\dif z\dif r\dif x\right|\\
&\qquad \preceq \int_{\mR^d}\Big(\varrho^0_{1+\beta}(t,x;s,y)+\varrho^{\beta}_{1}(t,x;s,y)\Big)(\dif y+\dif x)
\stackrel{(\ref{ES4})}{\preceq} |t-s|^{\beta}\rightarrow 0,\ t\downarrow s.
\end{align*}
Thus, (\ref{eq14}) and (\ref{eq141}) are proven by equation (\ref{eq2}).

For proving (\ref{EK1}), if we set $u_s(t,x):=\int_{\mR^d}p_{a,b}(t,x;s,y)\dif y$, then by (\ref{eq15}) and (\ref{eq14}),
$$
\p_tu_s(t,x)=\sL^x_{a,b}(t,x)u_s(t,x),\ \ \lim_{t\downarrow s}u_s(t,x)=1.
$$
By the maximal principal of nonlocal equation (cf. \cite{Zh1} or \cite[Theorem 2.3]{Zh3}), it follows that
$$
u_s(t,x)\equiv 1,\ \ t>s,\ \ x\in\mR^d.
$$

(iii) By (\ref{eq7}), one has
\begin{align}
\int_s^t\!\!\!\int_{\mR^d}p_0(t,x;r,z)|q(r,z;s,y)|\dif z\dif r
\preceq\varrho^0_{1+\beta}(t,x;s,y)+\varrho^{\beta}_{1}(t,x;s,y)\leq\varrho^0_1(t,x;s,y),\label{eq20}
\end{align}
which in turn gives estimate (\ref{eq16}) by equation (\ref{eq2}) and (\ref{E1}).

In the case of $a(t,x)=a(t)$, by (\ref{eq2}) and (\ref{eq5}), we have
\begin{align*}
|p(t,x;s,y)-p_0(t,x;s,y)|&\preceq \int_s^t\!\!\!\int_{\mR^d}p_0(t,x;r,z)q(r,z;s,y)\dif z\dif r\\
&\preceq \int_s^t\!\!\!\int_{\mR^d}\varrho^0_1(t,x;r,z)\varrho^0_{\beta}(r,z;s,y)\dif z\dif r\\
&\preceq |t-s|^{\beta}\varrho^0_1(t,x;s,y)\leq \Lambda|t-s|^{\beta}p_0(t,x;s,y),
\end{align*}
where $\Lambda>0$ is a constant independent of $t,s,x,y$.
Choosing $T_0\in(0,1)$ such that $|t-s|^{\beta}<\frac{1}{2\Lambda}$ for all $0\leq s<t\leq T_0$,
we obtain (\ref{eq18}) for small time by (\ref{E1}). For the large time, it follows by a standard time shift argument
(see \cite{Bo-Ja0, Wa-Zh}).

(iv) As in proving (\ref{ES5}), we have for any $\gamma\in(0,1)$,
\begin{align*}
|p_{0}(t,x;s,y)-p_{0}(t,x';s,y)|\preceq (|x-x'|^{\gamma}\wedge 1)
\Big(\varrho_{1-\gamma}^0(t,x;s,y)+\varrho_{1-\gamma}^0(t,x';s,y)\Big).
\end{align*}
Thus, by (\ref{eq3}) and Lemma \ref{Le2}, we have
\begin{align*}
&\int_s^t\!\!\!\int_{\mR^d}|p_{0}(t,x;r,z)-p_{0}(t,x';r,z)||q(r,z;s,y)|\dif z\dif r\\
&\quad\preceq\Big(|x-x'|^{\gamma}\wedge 1\Big)\int_s^t\!\!\!\int_{\mR^d}
\Big(\varrho_{1-\gamma}^0(t,x;r,z)+\varrho_{1-\gamma}^0(t,x';r,z)\Big)
(\varrho_{\beta}^0+\varrho^{\beta}_0)(r,z;s,y)\dif z\dif r\\
&\quad\preceq\Big(|x-x'|^{\gamma}\wedge 1\Big)
\Big((\varrho_{1+\beta-\gamma}^0+\varrho_{1-\gamma}^{\beta})(t,x;s,y)
+(\varrho_{1+\beta-\gamma}^0+\varrho_{1-\gamma}^{\beta})(t,x';s,y)\Big)\\
&\quad\preceq (|x-x'|^{\gamma}\wedge 1)\Big(\varrho_{1-\gamma}^0(t,x;s,y)+\varrho_{1-\gamma}^0(t,x';s,y)\Big),
\end{align*}
which together with equation (\ref{eq2}) yields (\ref{eq19}).

Next, we prove (\ref{eq17}). By (\ref{EP2}), we can write
\begin{align*}
\nabla_x \varphi(t,x)&=\int^t_{\frac{t+s}{2}}\!\int_{\mR^d}\nabla_x p_0(t,x;r,z)\Big(q(r,z;s,y)-q(r,x;s,y)\Big)\dif z\dif r\\
&\quad+\int^t_{\frac{t+s}{2}}\left(\int_{\mR^d}\nabla_x p_0(t,x;r,z)\dif z\right)q(r,x;s,y)\dif r\\
&\quad+\int^{\frac{t+s}{2}}_s\!\!\!\!\int_{\mR^d}\nabla_x p_0(t,x;r,z)q(r,z;s,y)\dif z\dif r\\
&=:Q_1(t,x;s,y)+Q_2(t,x;s,y)+Q_3(t,x;s,y).
\end{align*}
For $Q_1(t,x;s,y)$, by (\ref{EU4}), (\ref{eq4}) and Lemma \ref{Le2}, we have
\begin{align*}
|Q_1(t,x;s,y)|&\preceq\int^t_{\frac{t+s}{2}}\!\int_{\mR^d}
\varrho^{\beta-\gamma}_0(t,x;r,z)\Big\{(\varrho^0_\gamma+\varrho^\beta_{\gamma-\beta})(r,x;s,y)
+(\varrho^0_\gamma+\varrho^\beta_{\gamma-\beta})(r,z;s,y)\Big\}\dif z\dif r\\
&\preceq\int^t_{\frac{t+s}{2}}\left(\int_{\mR^d}\varrho^{\beta-\gamma}_0(t,x;r,z)\dif z\right)
(\varrho^0_\gamma+\varrho^\beta_{\gamma-\beta})(r,x;s,y)\dif r\\
&\quad+\int^t_s\!\!\!\int_{\mR^d}\varrho^{\beta-\gamma}_0(t,x;r,z)(\varrho^0_\gamma+\varrho^\beta_{\gamma-\beta})(r,z;s,y)\dif z\dif r\\
&\stackrel{(\ref{ES4})(\ref{eq30})}{\preceq}\left(\int^t_{\frac{t+s}{2}}(t-r)^{\beta-\gamma-1}
(1+(r-s)^{\gamma-\beta})\varrho^0_0(r,x;s,y)\dif r\right)\\
&\quad+(\varrho^0_\beta+\varrho^\beta_{0}+\varrho^{\beta-\gamma}_{\gamma})(t,x;s,y)\preceq\varrho^0_0(t,x;s,y).
\end{align*}
For $Q_2(t,x;s,y)$,  we have
\begin{align*}
|Q_2(t,x;s,y)|\stackrel{(\ref{EW7})(\ref{eq3})}{\preceq}\int^t_{\frac{t+s}{2}}(t-r)^{\beta-1}
\Big\{\varrho^\beta_0(r,x;s,y)+\varrho^0_\beta(r,x;s,y)\Big\}\dif r
\preceq\varrho^0_0(t,x;s,y).
\end{align*}
For $Q_3(t,x;s,y)$, we have
\begin{align*}
|Q_3(t,x;s,y)|\stackrel{(\ref{EU4})(\ref{eq3})}{\preceq}\int^{\frac{t+s}{2}}_s\!\!\!\!\int_{\mR^d}
\varrho^0_0(t,x;r,z)\Big\{\varrho^\beta_0(r,z;s,y)+\varrho^0_\beta(r,z;s,y)\Big\}\dif z\dif r
\stackrel{(\ref{EU7})(\ref{ES4})}{\preceq} \varrho^0_0(t,x;s,y).
\end{align*}
Combining the above calculations, we obtain
\begin{align}
|\nabla_x \varphi(t,x)|\preceq\varrho^0_0(t,x;s,y).\label{EQ5}
\end{align}
Similarly,
\begin{align}
|\Delta^{\frac{1}{2}}_x \varphi(t,x)|\preceq\varrho^0_0(t,x;s,y).\label{EQ6}
\end{align}
Estimate (\ref{eq17}) then follows by equation (\ref{eq2}), (\ref{EU5}), (\ref{EU4}) and (\ref{EQ5}), (\ref{EQ6}).

(v) For $f\in C^\infty_0(\mR^d)$, by (\ref{eq15}) and (\ref{eq14}), we have
\begin{align*}
\frac{P^{a,b}_{t,s}f(x)-f(x)}{t-s}-\sL^x_{a,b}(s,x)f(x)&=\frac{1}{t-s}\int^t_s(a(r,x)-a(s,x))
\Delta^{\frac{1}{2}}_xP^{a,b}_{r,s}f(x)\dif r\\
&\quad+\frac{a(s,x)}{t-s}\int^t_s\Delta^{\frac{1}{2}}_x(P^{a,b}_{r,s}f(x)-f(x))\dif r\\
&\quad+\frac{1}{t-s}\int^t_s(b(r,x)-b(s,x))\cdot\nabla_x P^{a,b}_{r,s}f(x)\dif r\\
&\quad+\frac{b(s,x)}{t-s}\cdot\int^t_s\nabla_x (P^{a,b}_{r,s}f(x)-f(x))\dif r\\
&=:I_1(t,s,x)+I_2(t,s,x)+I_3(t,s,x)+I_4(t,s,x).
\end{align*}
We have the following claim: for $t>s$,
\begin{align}
\sup_{r\in[s,t]}\|\Delta^{\frac{1}{2}}_xP^{a,b}_{r,s}f\|_\infty
+\sup_{r\in[s,t]}\|\nabla_xP^{a,b}_{r,s}f\|_\infty<+\infty,\label{EK2}
\end{align}
{\it Proof of Claim}: By Lemma \ref{Le6} and (\ref{EK1}), we have
$$
\Delta^{\frac{1}{2}}_xP^{a,b}_{r,s}f(x)=\int_{\mR^d}\Delta^{\frac{1}{2}}_xp_{a,b}(r,x;s,y)f(y)\dif y
=\int_{\mR^d}\Delta^{\frac{1}{2}}_xp_{a,b}(r,x;s,y)(f(y)-f(x))\dif y
$$
and
$$
\nabla_xP^{a,b}_{r,s}f(x)=\int_{\mR^d}\nabla_xp_{a,b}(r,x;s,y)f(y)\dif y
=\int_{\mR^d}\nabla_xp_{a,b}(r,x;s,y)(f(y)-f(x))\dif y.
$$
Hence, by (\ref{eq17}) we have
\begin{align*}
|\Delta^{\frac{1}{2}}_xP^{a,b}_{r,s}f(x)|
+|\nabla_xP^{a,b}_{r,s}f(x)|\leq C\|f\|_{\mathrm{Lip}}\int_{\mR^d}\varrho^1_0(r,x;s,y)\dif y
\stackrel{(\ref{ES4})}{\leq} C,
\end{align*}
which gives (\ref{EK2}).

\vspace{5mm}

For $g\in C^\infty_0(\mR^d)$, since $a,b\in C([0,\infty);L^1_{\mathrm{loc}}(\mR^d))$, by (\ref{EK2}), we have
\begin{align*}
\lim_{t\downarrow s}\int_{\mR^d}g(x)(I_1(t,s,x)+I_3(t,s,x))\dif x=0.
\end{align*}
Let $a_n(x)\in C^\infty_b(\mR^d)$ with
$$
\lim_{n\to\infty}|a_n(x)-a(s,x)|=0,\ \ \ x\in\mR^d.
$$
Then, by (\ref{EK2}), we have
\begin{align*}
\lim_{t\downarrow s}\left|\int_{\mR^d}g(x)I_2(t,s,x)\dif x\right|
&\leq\lim_{n\to\infty}\lim_{t\downarrow s}\left|\int_{\mR^d}g(x)\frac{a_n(x)-a(s,x)}
{t-s}\int^t_s\Delta^{\frac{1}{2}}_x(P^{a,b}_{r,s}f(x)-f(x))\dif r\dif x\right|\\
&+\lim_{n\to\infty}\lim_{t\downarrow s}\left|\int_{\mR^d}g(x)\frac{a_n(x)}
{t-s}\int^t_s\Delta^{\frac{1}{2}}_x(P^{a,b}_{r,s}f(x)-f(x))\dif r\dif x\right|\\
&\preceq\lim_{n\to\infty}\int_{\mR^d}|g(x)|\cdot|a_n(x)-a(s,x)|\dif x\\
&+\lim_{n\to\infty}\lim_{t\downarrow s}\left|\frac{1}
{t-s}\int^t_s\!\!\!\int_{\mR^d}\Delta^{\frac{1}{2}}_x(ga_n)(x)
(P^{a,b}_{r,s}f(x)-f(x))\dif x\dif r\right|,
\end{align*}
which converges to zero by (\ref{eq14}). Similarly, we also have
$$
\lim_{t\downarrow s}\left|\int_{\mR^d}g(x)I_4(t,s,x)\dif x\right|=0.
$$
Combining the above calculations, we obtain (\ref{eq233}).

(vi) Lastly, we show the uniqueness. For $f\in C^\infty_0(\mR^d)$ and $t\geq s$, set
$$
u^f_s(t,x):=\int_{\mR^d}p_{a,b}(t,x;s,y)f(y)\dif y.
$$
Let us first show the following claim: For any $T>s$ and $p>1$,
\begin{align}
\sup_{t\in[s,T]}\Big(\|u^f_s(t)\|_p+\|\nabla u^f_s(t)\|_p\Big)<+\infty.\label{ET1}
\end{align}
{\it Proof of Claim}: By Young's inequality for convolution, it follows that
$$
\|u^f_s(t)\|^p_p\stackrel{(\ref{eq16})}{\preceq}\int_{\mR^d}\left(\int_{\mR^d}
\varrho^0_1(t,x;s,y)|f(y)|\dif y\right)^p\dif x\stackrel{(\ref{ES4})}{\preceq}\|f\|^p_p.
$$
We next prove $\sup_{t\in[s,T]}\|\nabla u^f_s(t)\|_p<+\infty$.
Let the support of $f$ be contained in the ball $\{x\in\mR^d: |x|\leq N\}$.  We have
\begin{align*}
\int_{\mR^d}\left|\int_{\mR^d}\nabla_xp_{a,b}(t,x;s,y)f(y)\dif y\right|^p\dif x
&\preceq I_1+I_2,
\end{align*}
where
\begin{align*}
I_1&:=\int_{|x|\geq 2N}\left|\int_{\mR^d}\nabla_xp_{a,b}(t,x;s,y)f(y)\dif y\right|^p\dif x,\\
I_2&:=\int_{|x|< 2N}\left|\int_{\mR^d}\nabla_xp_{a,b}(t,x;s,y)f(y)\dif y\right|^p\dif x.
\end{align*}
For $I_1$, we have
\begin{align*}
I_1&=\int_{|x|\geq 2N}\left|\int_{|y|\leq N}\nabla_xp_{a,b}(t,x;s,y)f(y)\dif y\right|^p\dif x\\
&\leq C_N\int_{|x|\geq 2N}\int_{|y|\leq N}|\nabla_xp_{a,b}(t,x;s,y)|^p|f(y)|^p\dif y\dif x\\
&\leq C_N\|f\|^p_\infty\int_{|y|\leq N}\int_{|x-y|\geq N}|\nabla_xp_{a,b}(t,x;s,y)|^p\dif x\dif y\\
&\stackrel{(\ref{eq17})}{\leq}
C_N\|f\|^p_\infty\int_{|y|\leq N}\int_{|x-y|\geq N}\frac{1}{|x-y|^{(d+1)p}}\dif x\dif y<+\infty.
\end{align*}
For $I_2$, by (\ref{EK2}), we have
$$
I_2\leq\sup_{x\in\mR^d}\left|\int_{\mR^d}\nabla_xp_{a,b}(t,x;s,y)f(y)\dif y\right|^p\left(\int_{|x|<2N}\dif x\right)<+\infty.
$$
The claim is proven.

\vspace{5mm}

Let $p'_{a,b}(t,x;s,y)$ be another function satisfying (\ref{EK1})-(\ref{eq17}).
We want to prove that for any $f\in C^\infty_0(\mR^d)$ and $t>s$,
$$
\tilde u^f_s(t,x):=\int_{\mR^d}p_{a,b}(t,x;s,y)f(y)\dif y-\int_{\mR^d}p'_{a,b}(t,x;s,y)f(y)\dif y=0.
$$
In view of $\lim_{t\downarrow s}u^f_s(t,x)=0$, by (\ref{eq15}) we have
$$
\tilde u^f_s(t,x)=\int^t_s\sL^x_{a,b}(r,x)\tilde u^f_s(r,\cdot)(x)\dif r.
$$
The uniqueness follows by (\ref{ET1}) and \cite[Lemma 3.1]{Zh3}.
\end{proof}

\section{Proof of Theorem \ref{Main}}

By Duhamel's formula, we construct the heat kernel $p(t,x;s,y)$ of $\sL(t,x)$ by solving the following integral equation:
\begin{align}
p(t,x;s,y)=p_{a,b}(t,x;s,y)+\int_s^t\!\!\!\int_{\mR^d}p_{a,b}(t,x;r,z)c(r,z)p(r,z;s,y)\dif z\dif r. \label{eq8}
\end{align}
For $t>s\geq 0$ and $x,y\in\mR^d$, set $\Theta_0(t,x;s,y):=p_{a,b}(t,x;s,y)$, and define recursively for $n\in\mN$,
\begin{align}
\Theta_n(t,x;s,y)&:=\int^t_s\!\!\!\int_{\mR^d}p_{a,b}(t,x;r,z)c(r,z)\Theta_{n-1}(r,z;s,y)\dif z \dif r.\label{eq33}
\end{align}
For $\gamma\in(0,1]$ and $c\in\mK^\gamma_d$, define
$$
\ell^c_\gamma(\eps):=\sup_{(t,x)\in[0,\infty)\times\mR^d}\int^\eps_0\!\!\!\int_{\mR^d}
\varrho^0_\gamma(s,z)(|c(t-s,x-z)|+|c(t+s,x+z)|) \dif z\dif s.
$$
\bl\label{Le5}
If $c\in\mK^1_d$, then there exists a constant $\Lambda>0$ such that for all $n\in\mN$,
\begin{equation}\label{eq10}
|\Theta_n(t,x;s,y)|\leq\{\Lambda\ell^c_1 (t-s)\}^n\varrho^0_1(t,x;s,y).
\end{equation}
If $c\in \mK^{1-\gamma}_d$ for some $\gamma\in(0,1)$, then there exists a constant $C_1>0$
such that for any $n\in\mN$,
\begin{align}\label{eq35}
|\Theta_n(t,x;s,y)-\Theta_n(t,x';s,y)|&\leq C_1(|x-x'|^\gamma\wedge 1)
\{\Lambda\ell^c_1(t-s)\}^{n-1}\ell^c_{1-\gamma}(t-s)\no\\
&\qquad\times(\varrho^0_1(t,x;s,y)+\varrho^0_1(t,x';s,y)).
\end{align}
If $c\in\mH^{\gamma}$ for some $\gamma\in(0,1)$, then there exists a constant $C_2>0$
such that for any $n\in\mN$,
\begin{equation}\label{eq12}
|\nabla_x\Theta_n(t,x;s,y)|\leq C_2\{\Lambda\|c\|_\infty(t-s)\}^n\varrho_0^0(t,x;s,y).
\end{equation}
\el
\begin{proof}
(1) First of all, by (\ref{eq16}), we have for some $C_0>0$,
$$
p_{a,b}(t,x;s,y)\leq C_0\varrho^0_1(t,x;s,y).
$$
Now we use induction to prove (\ref{eq10}). Suppose that (\ref{eq10}) is true for $n\in\mN$.
Then
\begin{align*}
|\Theta_{n+1} (t,x;s,y)|&\leq \int^t_s\!\!\!\int_{\mR^d}p_{a,b}(t,x;r,z)|c(r,z)|\cdot |\Theta_{n} (r,z;s,y)|\dif z\dif r\\
&\leq C_0\{\Lambda\ell^c_1 (t-s)\}^n \int^t_s\!\!\!\int_{\mR^d}\varrho^0_1(t,x;r,z)\varrho^0_1(r,z;s,y)|c(r,z)|\dif z\dif r\\
&\stackrel{(\ref{eq6})}{\leq} \Lambda\{\Lambda\ell^c_1 (t-s)\}^n
\int^t_s\!\!\!\int_{\mR^d}\big(\varrho^0_1(t,x;r,z)+\varrho^0_1(r,z;s,y)\big) |c(r,z)| \dif z\dif r\varrho^0_1(t,x;s,y)\\
&\leq \{\Lambda\ell^c_1(t-s)\}^{n+1}\varrho^0_1(t,x;s,y).
\end{align*}

(2) By (\ref{eq33}) and (\ref{eq19}), we have
\begin{align*}
&|\Theta_n(t,x;s,y)-\Theta_n(t,x';s,y)|\\
&\quad\preceq(|x-x'|^{\gamma}\wedge 1)\int^t_s\!\!\!\int_{\mR^d}
\Big(\varrho_{1-\gamma}^0(t,x;r,z)+\varrho_{1-\gamma}^0(t,x';r,z)\Big)
|c(r,z)|\cdot|\Theta_{n-1}(r,z;s,y)|\dif z \dif r\\
&\quad\stackrel{(\ref{eq10})}{\preceq}(|x-x'|^{\gamma}\wedge 1)\{\Lambda\ell^c_1(t-s)\}^{n-1}\\
&\qquad\qquad\times\int^t_s\!\!\!\int_{\mR^d}\Big(\varrho_{1-\gamma}^0(t,x;r,z)+\varrho_{1-\gamma}^0(t,x';r,z)\Big)
|c(r,z)|\varrho^0_1(r,z;s,y)|\dif z \dif r\\
&\quad\stackrel{(\ref{ET2})}{\preceq}(|x-x'|^{\gamma}\wedge 1)\{\Lambda\ell^c_1(t-s)\}^{n-1}\\
&\qquad\qquad\times\Bigg\{\int^t_s\!\!\!\int_{\mR^d}(t-r)^{1-\gamma}(r-s)\Big(\varrho_0^0(t,x;r,z)+\varrho^0_0(r,z;s,y)\Big)
|c(r,z)|\dif z\dif r\rho^0_0(t,x;s,y)\\
&\qquad\qquad+\int^t_s\!\!\!\int_{\mR^d}(t-r)^{1-\gamma}(r-s)\Big(\varrho_0^0(t,x';r,z)+\varrho^0_0(r,z;s,y)\Big)
|c(r,z)|\dif z\dif r\rho^0_0(t,x';s,y)\Bigg\}\\
&\quad\leq C_1(|x-x'|^{\gamma}\wedge 1)\{\Lambda\ell^c_1(t-s)\}^{n-1}\\
&\qquad\qquad\times\Bigg\{\int^t_s\!\!\!\int_{\mR^d}\Big(\varrho_{1-\gamma}^0(t,x;r,z)+\varrho^0_{1-\gamma}(r,z;s,y)\Big)
|c(r,z)|\dif z\dif r\rho^0_1(t,x;s,y)\\
&\qquad\qquad+\int^t_s\!\!\!\int_{\mR^d}\Big(\varrho_{1-\gamma}^0(t,x';r,z)+\varrho^0_{1-\gamma}(r,z;s,y)\Big)
|c(r,z)|\dif z\dif r\rho^0_1(t,x';s,y)\Bigg\}\\
&\quad\leq C_1(|x-x'|^{\gamma}\wedge 1)\{\Lambda\ell^c_1(t-s)\}^{n-1}\ell^c_{1-\gamma}(t-s)(\varrho_1^0(t,x;s,y)+\varrho_1^0(t,x';s,y)),
\end{align*}
and (\ref{eq35}) holds.

(3) If $c$ is bounded, by definition and (\ref{ES4}), it is easy to see that for some $C_1>0$,
\begin{align}
\ell^c_{\gamma}(\eps)\leq C_1\|c\|_\infty\eps^\gamma, \ \eps>0.\label{ET3}
\end{align}
As in Lemma \ref{Le6}, one can prove
$$
\nabla_x\Theta_n(t,x;s,y)=\int^t_s\!\!\!\int_{\mR^d}\nabla_xp_{a,b}(t,x;r,z)c(r,z)\Theta_{n-1}(r,z;s,y)\dif z \dif r.
$$
By (\ref{EK1}), we can write
\begin{align*}
\nabla_x\Theta_n(t,x;s,y)&=\int^t_{\frac{t+s}{2}}\!\int_{\mR^d}\nabla_xp_{a,b}(t,x;r,z)\Big(c(r,z)\Theta_{n-1}(r,z;s,y)
-c(r,x)\Theta_{n-1}(r,x;s,y)\Big)\dif z \dif r\\
&\quad+\int^{\frac{t+s}{2}}_s\!\!\!\int_{\mR^d}\nabla_xp_{a,b}(t,x;r,z)c(r,z)\Theta_{n-1}(r,z;s,y)\dif z \dif r\\
&=\int^t_{\frac{t+s}{2}}\!\int_{\mR^d}\nabla_xp_{a,b}(t,x;r,z)c(r,z)\Big(\Theta_{n-1}(r,z;s,y)
-\Theta_{n-1}(r,x;s,y)\Big)\dif z\dif r\\
&\quad+\int^t_{\frac{t+s}{2}}\left(\int_{\mR^d}\nabla_xp_{a,b}(t,x;r,z)(c(r,z)-c(r,x))\dif z\right)
\Theta_{n-1}(r,x;s,y)\dif r\\
&\quad+\int^{\frac{t+s}{2}}_s\!\!\!\int_{\mR^d}\nabla_xp_{a,b}(t,x;r,z)c(r,z)\Theta_{n-1}(r,z;s,y)\dif z \dif r\\
&=:Q_1(t,x;s,y)+Q_2(t,x;s,y)+Q_3(t,x;s,y).
\end{align*}
For $Q_1(t,x;s,y)$, by (\ref{ET3}) and (\ref{eq35}), we have
\begin{align*}
Q_1(t,x;s,y)&\preceq
\{\Lambda\|c\|_\infty(t-s)\}^{n-1}\int_{\frac{t+s}{2}}^t\!\int_{\mR^d}\varrho^{\gamma}_0(t,x;r,z)\varrho^0_{1-\gamma}(r,z;s,y)\dif z\dif r\\
&\quad +\{\Lambda\|c\|_\infty(t-s)\}^{n-1}\int_{\frac{t+s}{2}}^t\left(\int_{\mR^d}\varrho^{\gamma}_0(t,x;r,z)\dif z\right)
\varrho^0_{1-\gamma}(r,x;s,y)\dif r\\
&\stackrel{(\ref{EU7})(\ref{ES4})}{\preceq}\{\Lambda\|c\|_\infty(t-s)\}^n\varrho_0^0(t,x;s,y).
\end{align*}
For $Q_2(t,x;s,y)$, by (\ref{ET3}) and (\ref{eq10}), we have
\begin{align*}
Q_2(t,x;s,y)&\preceq
\{\Lambda\|c\|_\infty(t-s)\}^{n-1}\int^t_{\frac{t+s}{2}}\left(\int_{\mR^d}\varrho^{\gamma}_0(t,x;r,z)\dif z
\right)\varrho^0_1(r,x;s,y)\dif r\\
&\preceq\{\Lambda\|c\|_\infty(t-s)\}^{n-1}\left(\int^t_{\frac{t+s}{2}}(t-r)^{\gamma-1}(r-s)\dif r\right)\varrho^0_0(t,x;s,y)\\
&\preceq\{\Lambda\|c\|_\infty(t-s)\}^{n}\varrho_0^0(t,x;s,y).
\end{align*}
For $Q_3(t,x;s,y)$, we have
\begin{align*}
Q_3(t,x;s,y)&\preceq
\{\Lambda\|c\|_\infty(t-s)\}^{n-1}\int^{\frac{t+s}{2}}_s\!\!\!\int_{\mR^d}\varrho^0_0(t,x;r,z)\varrho^0_1(r,z;s,y)\dif z \dif r\\
&\preceq\{\Lambda\|c\|_\infty(t-s)\}^{n-1}\left(\int^{\frac{t+s}{2}}_s((r-s)(t-r)^{-1}+1)\dif r\right)\varrho^0_0(t,x;s,y)\\
&\preceq\{\Lambda\|c\|_\infty(t-s)\}^{n}\varrho_0^0(t,x;s,y).
\end{align*}
Combining the above calculations, we obtain (\ref{eq12}).
\end{proof}

Now we are in a position to give

\begin{proof}[Proof of Theorem\ref{Main}] By the standard time shift technique, it suffices
to prove the conclusions on a small time interval. We divide the proof in several steps.

(1) Define
$$
p(t,x;s,y)=p_{a,b}(t,x;s,y)+\sum_{n=1}^\infty\Theta_n(t,x;s,y).
$$
By virtue of $c\in\mK^1_d$, we have
$$
\lim_{\eps\downarrow 0}\ell^c_1(\eps)=0.
$$
Hence, for any given $\eps\in(0,1)$, one can choose $T_\eps\in(0,1)$ small enough such that for all $0\leq s<t\leq T_\eps$,
$$
\ell^c_1(t-s)\leq\frac{\eps}{\Lambda}.
$$
Thus,
\begin{align*}
|p(t,x;s,y)-p_{a,b}(t,x;s,y)|\leq\sum_{n=1}^\infty|\Theta_n(t,x;s,y)|&\leq
\frac{\Lambda\ell^c_1(t-s)}{1-\Lambda\ell^c_1(t-s)}\varrho^0_1(t,x;s,y)\\
&\leq\frac{\eps}{1-\eps}\varrho^0_1(t,x;s,y),
\end{align*}
which together with (\ref{eq16}) gives (\ref{eq24}) for $0\leq s<t\leq T_\eps$. Moreover, noticing that
$$
\sum_{n=0}^m\Theta_n(t,x;s,y)=p_{a,b}(t,x;s,y)+\int^t_s\!\!\!\int_{\mR^d}p_{a,b}(t,x;r,z)c(r,z)\sum_{n=0}^{m-1}\Theta_n(r,z;s,y)\dif z \dif r,
$$
by taking limits, we obtain equation (\ref{eq8}).
In the case of $a(t,x)=a(t)$, by (\ref{eq18}), if we let $\eps$ be small enough, we also have (\ref{eq25}).
Moreover, estimates (\ref{eq34}) and (\ref{eq26}) follow by (\ref{eq35}), (\ref{eq19}) and
(\ref{eq12}), (\ref{eq17}).

(2) Define
$$
P_{t,s}f(x):=\int_{\mR^d}p(t,x;s,y)f(y)\dif y
$$
and
$$
P^{a,b}_{t,s}f(x):=\int_{\mR^d}p_{a,b}(t,x;s,y)f(y)\dif y.
$$
For proving (\ref{eq21}), it suffices to prove that for any $f\in C^\infty_0(\mR^d)$,
\begin{align}
P_{t,s}f(x)=P_{t,r}P_{r,s}f(x),\ \ s<r<t.\label{eq27}
\end{align}
By (\ref{eq8}), we have
\begin{align*}
P_{t,s}f(x)&=P^{a,b}_{t,s}f(x)+\int^t_s P^{a,b}_{t,r'}\Big(c(r',\cdot)P_{r',s}f\Big)(x)\dif r'\\
&=P^{a,b}_{t,r}P^{a,b}_{r,s}f(x)+\int_s^rP^{a,b}_{t,r}P^{a,b}_{r,r'}\Big(c(r',\cdot)P_{r',s}f\Big)(x)\dif r'
+\int^t_rP^{a,b}_{t,r'}\Big(c(r',\cdot)P_{r',s}f\Big)(x)\dif r'\\
&=P^{a,b}_{t,r}P_{r,s}f(x)+\int^t_rP^{a,b}_{t,r'}\Big(c(r',\cdot)P_{r',s}f\Big)(x)\dif r',
\end{align*}
where we have used $P^{a,b}_{t,s}f=P^{a,b}_{t,r}P^{a,b}_{r,s}f$, which follows by the uniqueness of Theorem \ref{T1}.
On the other hand, we also have
\begin{align*}
P_{t,r}P_{r,s}f(x)=P^{a,b}_{t,r}P_{r,s}f(x)+\int^t_rP^{a,b}_{t,r'}\Big(c(r',\cdot)P_{r',r}P_{r,s}f\Big)(x)\dif r'.
\end{align*}
Fix $s<r$ and set
$$
u_t(x):=P_{t,r}P_{r,s}f(x)-P_{t,s}f(x).
$$
Then, we have
$$
u_t(x)=\int^t_r\!\!\!\int_{\mR^d}p_{a,b}(t,x;r',y)c(r',y)u_{r'}(y)\dif y\dif r'.
$$
By (\ref{eq16}), we have
$$
\|u_t\|_\infty\leq\sup_{r'\in[r,t]}\|u_{r'}\|_\infty
\int^t_r\!\!\!\int_{\mR^d}\varrho^0_1(t,x;r',y)|c(r',y)|\dif y\dif r'=\ell^c_1(t-r)\sup_{r'\in[r,t]}\|u_{r'}\|_\infty,
$$
which implies that
$$
\sup_{r'\in[r,t]}\|u_{r'}\|_\infty\leq\sup_{\eps\in(0,t-r]}\ell^c_1(\eps)\sup_{r'\in[r,t]}\|u_{r'}\|_\infty.
$$
In particular, if $t-r$ is small enough (say less than $\eps_0$), then
$$
\sup_{r'\in[r,t]}\|u_{r'}\|_\infty=0.
$$
Thus, we obtain (\ref{eq27}) for $t-r<\eps_0$. For general $t$, it follows by repeatedly using (\ref{eq27}).

(3) We prove (\ref{eq22}). By (\ref{eq8}) and (\ref{eq14}), we only need to prove that for any $f\in C_b(\mR^d)$,
$$
\lim_{t\downarrow s}\int_{\mR^d}\!\int^t_s\!\!\!\int_{\mR^d}p_{a,b}(t,x;r,z)c(r,z)p(r,z;s,y)f(y)\dif z\dif r\dif y = 0.
$$
This limit follows by noticing that
\begin{align*}
&\left|\int_{\mR^d}\!\int^t_s\!\!\!\int_{\mR^d}p_{a,b}(t,x;r,z)c(r,z)p(r,z;s,y)f(y)\dif z\dif r\dif y\right|\\
& \quad \preceq \int_{\mR^d}\!\int^t_s\!\!\!\int_{\mR^d}\varrho^0_1(t,x;r,z)|c(r,z)|\varrho^0_1(r,z;s,y)|f(y)|\dif z\dif r\dif y\\
& \quad \preceq \int_{\mR^d}\left(\int^t_s\!\!\!\int_{\mR^d}\big(\varrho^0_1(t,x;r,z)+\varrho^0_1(r,z;s,y)\big)
c(r,z)\dif z\dif r\right)\varrho^0_1(t,x;s,y)\dif y\\
& \quad \preceq \ell^c_1(t-s)\int_{\mR^d}\varrho^0_1(t,x;s,y)\dif y\stackrel{(\ref{ES4})}{\leq}
C\ell^c_1(t-s)\to 0, \ t\downarrow s.
\end{align*}

(4) Let $f,g\in C^\infty_0(\mR^d)$. By definitions, we make the following decomposition:
\begin{align*}
\frac{P_{t,s}f(x)-f(x)}{t-s}-\sL f(x)&=\frac{1}{t-s}\int_s^t\Big(P^{a,b}_{t,r}\big(c(r)P_{r,s}f\big)(x)-c(r,x)P_{r,s}f(x)\Big)\dif r\\
&\quad+ \frac{1}{t-s}\int_s^t\Big(c(r,x)-c(s,x)\Big)P_{r,s}f(x)\dif r\\
&\quad+ \frac{1}{t-s}\int_s^tc(s,x)\Big(P_{r,s}f(x)-f(x)\Big)\dif r\\
&\quad+\Bigg(\frac{P^{a,b}_{t,s}f(x)-f(x)}{t-s}-\sL^x_{a,b}(s,x)f(x)\Bigg)\\
&=: I_1(t,s,x)+I_2(t,s,x)+I_3(t,s,x)+I_4(t,s,x).
\end{align*}
For $I_1(t,s,x)$, if we write
$$
(P^{a,b}_{t,r})^*g(y):=\int_{\mR^d}p_{a,b}(t,x;r,y)g(x)\dif x,
$$
then
\begin{align*}
\left|\int_{\mR^d}g(x)I_1(t,s,x)\dif x\right|&\leq\left|\frac{1}{t-s}\int_s^t\!\!\!\int_{\mR^d}
\Big((P^{a,b}_{t,r})^*g(x)-(P^{a,b}_{t,r})^*1(x)\cdot g(x)\Big)c(r,x)P_{r,s}f(x)\dif x\dif r\right|\\
&\quad+\left|\frac{1}{t-s}\int_s^t\!\!\!\int_{\mR^d}\Big((P^{a,b}_{t,r})^*1-1\Big)(x) g(x)c(r,x)P_{r,s}f(x)\dif x\dif r\right|\\
&=:J_1(t,s)+J_2(t,s).
\end{align*}
For $J_1(t,s)$, noticing that
\begin{align*}
|(P^{a,b}_{t,r})^*g(y)-(P^{a,b}_{t,r})^*1(y)\cdot g(y)|&=\left|\int_{\mR^d}p_{a,b}(t,x;r,y)(g(x)-g(y))\dif x\right|\\
&\stackrel{(\ref{eq16})}{\leq} C\|g\|_{\mH^1}\int_{\mR^d}\varrho^0_1(t,x;r,y)(|x-y|\wedge 1)\dif x\\
&\stackrel{(\ref{ES4})}{\leq} C\|g\|_{\mH^1}|t-r|,
\end{align*}
by definition of $P_{r,s}f$ and (\ref{eq24}), we have
\begin{align*}
J_1(t,s)&\leq C\|g\|_{\mH^1}\int_s^t\!\!\!\int_{\mR^d}|c(r,x)|\cdot|P_{r,s}f(x)|\dif x\dif r\\
&\leq C\|g\|_{\mH^1}\int_s^t\!\!\!\int_{\mR^d}\int_{\mR^d}|c(r,x)|\varrho^0_1(r,x;s,y)|f(y)|\dif y\dif x\dif r\\
&\leq C\|g\|_{\mH^1}\ell^c_1(t-s)\int_{\mR^d}|f(y)|\dif y\to 0,\ \ t\downarrow s.
\end{align*}
For $J_2(t,s)$, since $c\in C([0,\infty);L_{\mathrm{loc}}^{1}(\mR^d))$,
by (\ref{eq141}) and the dominated convergence theorem, we have
$$
\lim_{t\downarrow s}J_2(t,s)=0.
$$
It is the same reason that
$$
\lim_{t\downarrow s}\int_{\mR^d}g(x)(I_2(t,s,x)+I_3(t,s,x))\dif x=0.
$$
Moreover, if $a,b\in C([0,\infty);L_{\mathrm{loc}}^{1}(\mR^d))$, by (\ref{eq233}) we have
\begin{align*}
\lim_{t\downarrow s}\int_{\mR^d}g(x)I_4(t,s,x)\dif x=0.
\end{align*}
Combining the above limits, we obtain (\ref{eq23}). The whole proof is complete.
\end{proof}

{\bf Acknowledgements:}

The authors would like to thank Professors Zhen-Qing Chen, Renming Song and Feng-Yu Wang for
their quite useful conversations. This work is supported by NSFs of China (No. 11271294) and
Program for New Century Excellent Talents in University (NCET-10-0654).

\end{document}